\documentclass[preprint,12pt]{elsarticle}

\usepackage{graphicx}

\usepackage{amsmath}
\usepackage{amssymb}
\usepackage{amsthm}

\newtheorem{lem}{Lemma}
\newtheorem{prop}{Proposition}
\newtheorem{definition}{Definition}

\newtheorem{example}{Example}

\usepackage{array}
\usepackage{tikz}
\usepackage{pgfplots}
\usetikzlibrary{plotmarks,intersections,positioning,calc,trees}
\usepackage[latin1]{inputenc}
\usepackage{verbatim}
\usepackage[ruled, vlined]{algorithm2e}
\usepackage{caption}

\journal{Artificial Intelligence}

\begin{document}

\begin{frontmatter}



\title{Confidence-based Reasoning in \\Stochastic Constraint Programming\tnoteref{title_note}}

\tnotetext[title_note]{This work is an extended version of \cite{citeulike:9453442}}


\author[roberto_rossi]{Roberto Rossi\corref{author_note}}
\ead{roberto.rossi@ed.ac.uk}
\author[brahim_hnich] {Brahim Hnich}
\ead{hnich.brahim@gmail.com}
\author[armagan_tarim,steven_prestwich] {S. Armagan Tarim\fnref{steve_note}}
\ead{armtar@yahoo.com}
\author[steven_prestwich] {Steven Prestwich\fnref{steve_note}}
\ead{s.prestwich@cs.ucc.ie}

\address[roberto_rossi]{Business School, University of Edinburgh, United Kingdom}
\address[brahim_hnich]{Department of Computer Science, Taif University, Taif, Kingdom of Saudi Arabia}
\address[armagan_tarim]{Department of Management, Cankaya University, Ankara, Turkey}
\address[steven_prestwich]{Insight Centre for Data Analytics, University College Cork, Ireland}

\cortext[author_note]{{\bf Corresponding author}. Roberto Rossi,
Business School, University of Edinburgh, 
29 Buccleuch place, EH8 9JS, Edinburgh, UK.
Tel. +44(0)131 6515239, Fax. +44 (0)131 650 8077}

\fntext[steve_note]{
This publication has emanated from research supported in part by a research grant from Science Foundation Ireland (SFI) under Grant Number SFI/12/RC/2289}

\begin{abstract}
In this work we introduce a novel approach, based on sampling, for finding assignments that are likely to be solutions to stochastic constraint satisfaction problems and constraint optimisation problems. Our approach reduces the size of the original problem being analysed; by solving this reduced problem, with a given confidence probability, we obtain assignments that satisfy the chance constraints in the original model within prescribed error tolerance thresholds.  To achieve this, we blend concepts from stochastic constraint programming and statistics. We discuss both exact and approximate variants of our method. The framework we introduce can be immediately employed in concert with existing approaches for solving stochastic constraint programs. A thorough computational study on a number of stochastic combinatorial optimisation problems demonstrates the effectiveness of our approach.
\end{abstract}


\end{frontmatter}


\section{Introduction}\label{sec:int}

The stochastic constraint satisfaction/optimisation framework introduced in \cite{DBLP:conf/ecai/Walsh02,tmw2006} constitutes an expressive declarative formalism for modeling problems of decision making under uncertainty. A stochastic constraint satisfaction problem (SCSP), alongside decision variables, features \emph{random variables}, which 
follow some probability distribution and can be used to model uncertainty. Relationships over subsets of random and decision variables can be expressed in a declarative manner via \emph{stochastic constraints}. The fact that a given relationship over subsets of random and decision variables should be satisfied according to a prescribed probability can be expressed by means of \emph{chance constraints}. Finally, since problems of decision making under uncertainly are sequential in nature, the modeler can define a \emph{stage structure}, that is a sequence of \emph{decision stages}, in each of which a subset of all possible decisions are taken and a subset of all possible random variables are observed. 
A solution to an SCSP can be represented in general by means of a \emph{policy tree}, which records feasible or optimal decisions associated with each possible set of random variable realisations.

As shown in \cite[][Theorem 1]{tmw2006}, solving SCSPs is a computationally hard task. Even trivial instances with a dozen of decision and random variables require a computational effort out of reach even for the most advanced hardware/software combination. This is due to the fact that the size of the policy tree grows exponentially in the number of random variables in the model and in the size of their support. Furthermore, a major limitation of all existing complete SCSPs solution methods, such as \cite{tmw2006} and \cite{citeulike:10689458}, is the fact that they assume the support of random variables is \emph{finite}, otherwise a solution cannot be expressed as a finite policy tree. In practice, however, it is often the case that random variables either range over continuous supports or have a very large number (possibly infinite) of values in their domain. To date, no general purpose method exists for solving large-scale SCSPs, or SCSPs featuring random variables with continuous or discrete infinite support; for the sake of brevity we shall name this latter class of SCSPs ``infinite SCSPs.''

The main contribution of this paper is to propose a framework for solving large-scale or infinite SCSPs. More specifically, we argue that in solving large-scale or infinite SCSPs, one should not consider the ultimate feasible/optimal solution, which in some cases may even be impossible to represent; rather, the decision maker should aim for a solution that she ``sufficiently trusts,'' which she may claim to be optimal or feasible with a given confidence level, and for which a certain degree of error may be tolerated.  In order to obtain such a solution, the decision maker should only look at a possibly limited number of samples drawn from the random variables in the model. In other words, she should try to ``estimate'' the quality of this solution.

Our approach has several analogies with established techniques in statistics. 
When a survey is conducted on a sample population --- e.g. an electoral poll --- a statistician typically
associates a certain confidence level with the results obtained from the
chosen sample population. For instance, one may claim that
there is a 90\% chance that the actual mean being estimated is within a given
interval. We argue that the very same approach may be adopted in stochastic decision making.  
If the infinite or large-scale $m$-stage SCSP does not admit any closed form solution and is complex enough to rule out any chance of obtaining an exact solution, we suggest that --- as is done in statistics ---
one may introduce a confidence level $\alpha$ and a tolerated estimation error $\pm\vartheta$. 
The decision maker, instead of looking for an exact solution, 
may then aim to ``estimate'' --- according to the chosen $\alpha$ and $\vartheta$ --- whether the actual satisfaction probability guaranteed by an assignment is greater than or equal to the given target value for each of the chance constraints
in the model. 
By choosing given values for $\alpha$ and $\vartheta$ 
the set of solutions may vary. For this reason we will introduce a new notion of solution that is 
parameterised by these two parameters and that we call an $(\alpha,\vartheta)$-solution.
Intuitively, as $\alpha$ tends to 1 and $\vartheta$ tends to 0 the set of $(\alpha,\vartheta)$-solutions
will converge to the set of actual solutions to the original stochastic constraint satisfaction problem, 
which we therefore rename $(1,0)$-solutions. One should note that an approach of this kind
has been recently advocated in \cite[][chap. 4]{citeulike:13384427}.

In this work, we make the following contributions to the stochastic constraint programming literature:
\begin{itemize}
\item we discuss how to obtain compact instances of infinite or large-scale stochastic constraint programs via sampling: we call these instances ``sampled SCSPs;''
\item we introduce the concepts of $(\alpha,\vartheta)$-solution and of $(\alpha,\vartheta)$-solution set; and show how to compute a priori the minimum sample size that guarantees the attainment of such classes of solutions;
\item we show how the above tools can be employed in order to find approximate solutions to infinite or large-scale stochastic constraint satisfaction/optimisation problems that cannot be solved by existing exact approaches in the stochastic constraint programming literature;
\item we conduct a thorough computational study on three well-known stochastic combinatorial problems to validate our theoretical framework and assess its effectiveness, efficiency, and scalability.
\end{itemize}

This work is structured as follows: in Section \ref{sec:fb} we introduce the relevant formal background in constraint programming, stochastic constraint programming, and confidence interval analysis; in Section \ref{sec:sampled_scsps} we introduce sampled SCSPs; in Section \ref{sec:sampled_scsps_solutions} we discuss properties of the solutions of sampled SCSPs and formally introduce $(\alpha,\vartheta)$-solutions; in Section \ref{sec:alpha_theta_solution_set} we introduce $(\alpha,\vartheta)$-solution sets; in Section \ref{sec:scop} we extend our discussion to stochastic constraint optimisation problems; in Section \ref{sec:connections} we discuss connections with established techniques in statistics; in Section \ref{sec:comp} we present our computational study; in Sections \ref{sec:rel_works} we discuss related works; finally, in Section \ref{sec:con} we draw conclusions and discuss future research directions.

\section{Formal background}\label{sec:fb}

We now introduce the relevant background in  constraint programming, stochastic constraint programming, and confidence interval analysis.

\subsection{Constraint Programming}

A Constraint Satisfaction Problem (CSP) \cite{1207782} consists of a set of
decision variables, each with a finite domain of values, and a set of
constraints specifying allowed combinations of values for some
variables. A {\em solution} to a CSP is an assignment of variables
to values in their respective domains such that all of the
constraints are satisfied. Constraint solvers typically explore
partial assignments enforcing a local consistency property. A
constraint $c$ is \emph{generalized arc  consistent} (\emph{GAC})
if and only if when a variable is assigned any of the values in its domain,
there exist compatible values in the domains of all the other
variables of $c$. In order to enforce a local consistency
property on a constraint $c$ during search, we employ  filtering
algorithms that remove inconsistent values from the domains of the
variables of $c$. These filtering algorithms are repeatedly called
until no more values are pruned. This process is called
{\em constraint propagation}.

\subsection{Stochastic Constraint Programming}

The following definitions are based on \cite{DBLP:conf/cp/HnichRTP09,citeulike:10689458}. 
An $m$-stage stochastic constraint satisfaction problem (SCSP) \cite{DBLP:conf/ecai/Walsh02} 
is a 7-tuple $\langle V,S,D,P,C,\beta, L \rangle$, where $V$ is a set of decision
variables and $S$ is a set of random variables, $D$ is a
function mapping each element of $V$ (respectively, $S$) to a
domain (respectively, support) of potential values. In classical 
SCSPs both decision variable domains and random
variable supports are assumed to be finite. 
$P$ is a function mapping each element of $S$ to a
probability distribution for its associated support. To keep the discussion focused, 
without loss of generality, we assume that this probability distribution is not influenced by
the decisions made; extensions to the SCP framework that deal with decision-dependent probabilities
are discussed in \cite{citeulike:13377451}. $C$ is a set
of constraints over a non-empty subset of decision
variables and a subset of random variables. If a constraint
constrains only decision variables, then we call it
a deterministic constraint; if it constrains
both decision and random variables, then
we call it a stochastic constraint. $\beta$ is a
function mapping each stochastic constraint $h \in C$ to $\beta_h$,
which is a threshold value in the interval $(0,1]$. If this threshold
is strictly less than 1, then the stochastic constraint is a chance constraint.
$L=[\langle V_1,S_1
\rangle,\ldots,\langle V_i,S_i \rangle,\ldots, \langle V_m,S_m
\rangle]$ is  a list of {\em decision stages} such that each $V_i
\subseteq V$, each $S_i \subseteq S$, the $V_i$ form a partition
of $V$, and the $S_i$ form a partition of $S$. 

To solve an $m$-stage SCSP an assignment to the variables in $V_1$
must be found such that, given random values for $S_1$, assignments
can be found for $V_2$ such that, given random values for $S_2$,
$\ldots$, assignments can be found for $V_m$ so that, given random
values for $S_m$, the deterministic constraints are satisfied and the 
stochastic constraints are satisfied in the fraction of all possible
scenarios  specified by function $\beta$. Under the assumption that random variable
supports are finite, the solution of an $m$-stage SCSP is, in general, represented by
means of a {\em policy tree} \cite{tmw2006}. The arcs in such a
policy tree represent values observed for random variables
whereas  nodes at each level  represent the decisions associated
with the different stages. We call the policy tree of an
$m$-stage SCSP that is a solution a {\em satisfying policy tree}.

Let $\mathcal{S}$ denote the space of policy trees that are 
solutions to an SCSP. We may  be interested in
finding a policy tree $s\in \mathcal{S}$ that maximizes the value of a given objective
function $f(\cdot)$ over a subset of stochastic variables and a 
non-empty subset of decision variables. A 
stochastic constraint optimization problem (SCOP) is then 
defined in general as $\max_{s\in \mathcal{S}} f(s)$.

In order to simplify the presentation, we assume without loss of
generality, that each $V_i=\{x_i\}$ and each $S_i=\{s_i\}$ are
singleton sets. All the results can be easily 
extended in order to consider $|V_i| > 1$ and $|S_i|> 1$ (see \cite{citeulike:10689458}).
Intuitively, if $S_i$ comprises more than one random variable, it is always possible
to aggregate these variables into a single multivariate random variable \cite{j61} 
by convolving them. If $V_i$ comprises more than one decision variable, 
in the following discussion the term decision variable should be interpreted as a set of decision variables. 
Let $S=\{s_1,s_2,\ldots,s_m\}$ be the set of all random
variables and $V=\{x_1,x_2,\ldots,x_m\}$ be the set of all decision
variables. 

Let $p$ be a path from the root node of the policy tree
to a leaf. Let $\Psi$ denote the set of all distinct paths of a policy tree.
For each $p \in \Psi$, we denote by $\text{arcs}(p)$ the sequence of
all the arcs in $p$ whereas $\text{nodes}(p)$ denotes the
sequence of all nodes in $p$. We denote by $\Omega=\{\text{arcs}(p) | p \in
\Psi\}$ the set of all scenarios of the policy tree. The probability
of $\omega \in \Omega$ is given by $\Pr\{\omega\}=\prod_{i=1}^m
\Pr\{s_i=\bar{s}_i|s_{i-1}=\bar{s}_{i-1},\ldots,s_{1}=\bar{s}_{1}\}$, 
where $\Pr\{s_i=\bar{s}_i|s_{i-1}=\bar{s}_{i-1},\ldots,s_{1}=\bar{s}_{1}\}$ is the probability that 
random variable $s_i$ takes value $\bar{s}_i$, given a set of realisations for random variables 
$s_{i-1},\ldots,s_{1}$ already observed.

Now consider a constraint $h \in C$ with a specified
threshold level $\beta_h$. Consider a policy tree $\mathcal{T}$
for the SCSP and a path $p \in \mathcal{T}$. Let $h_{\downarrow{p}}$ 
be the deterministic constraint obtained by substituting the random 
variables in $h$ with the corresponding values ($\bar{s}_i$) assigned to 
these random variables in $\text{arcs}(p)$. Let $\bar{h}_{\downarrow{p}}$ 
be the resulting tuple obtained by substituting the decision variables 
in $h_{\downarrow{p}}$ by the values ($\bar{x}_i$) assigned to the 
corresponding decision variables in $\text{nodes}(p)$.
We say that $h$ is {\em satisfied wrt to a given policy tree $\mathcal{T}$} iff
\[\sum_{p \in \Psi: \bar{h}_{\downarrow{p}} \in h_{\downarrow{p}}} \Pr\{\text{arcs}(p)\} \geq
\beta_h.\]

\begin{definition}
Given an $m$-stage SCSP $\mathcal{P}$ and a policy tree
$\mathcal{T}$, $\mathcal{T}$ is a satisfying policy tree to
$\mathcal{P}$  iff every constraint of $\mathcal{P}$ is
satisfied wrt $\mathcal{T}$.
\end{definition}

\begin{figure}[t!]
\centering 
\tikzstyle{level 1}=[level distance=3.5cm, sibling distance=3.5cm]
\tikzstyle{level 2}=[level distance=3.5cm, sibling distance=2cm]

\tikzstyle{bag} = [text width=4em, text centered]
\tikzstyle{end} = [circle, minimum width=3pt,fill, inner sep=0pt]


\begin{tikzpicture}[grow=right]
\draw[name path=V1a--V1b,style=thin,color=gray] (1,-3) -- (1,3);
\draw[name path=V1a--V1b,style=thin,color=gray] (2.7,-3) -- (2.7,3);
\draw[name path=V1a--V1b,style=thin,color=gray] (4.7,-3) -- (4.7,3);
\draw[name path=V1a--V1b,style=thin,color=gray] (6.5,-3) -- (6.5,3);
\draw[name path=V1a--V1b,style=thin,color=gray] (8.3,-3) -- (8.3,3);
\draw (0,3.5) node[above] {\large$V_1$};
\draw (1.8,3.5) node[above] {\large$S_1$};
\draw (3.5,3.5) node[above] {\large$V_2$};
\draw (5.5,3.5) node[above] {\large$S_2$};
\draw (7.3,3.5) node[above] {Scenario prob.};
\node[bag] {$x_1=3$}
    child {
        node[bag] {$x_2^1=4$}        
            child {
                node[end, label=right:
                    {
                    0.25~~~~~
                    $\begin{array}{l}
                    c_1:5\cdot 3+4\cdot 4\geq 30\\
                    c_2:4\cdot 3 = 12
                    \end{array}$
                    }] {}
                edge from parent
                node[above] {$s_2=4$}
                node[below]  {$0.5$}
            }
            child {
                node[end, label=right:
                    {
                    0.25~~~~~
                    $\begin{array}{l}
                    c_1:5\cdot 3+3\cdot 4< 30\\
                    c_2:3\cdot 3 \neq 12
                    \end{array}$
                    }] {}
                edge from parent
                node[above] {$s_2=3$}
                node[below]  {$0.5$}
            }
            edge from parent 
            node[above] {$s_1=5$}
            node[below]  {$0.5$}
    }
    child {
        node[bag] {$x_2^2=6$}        
        child {
                node[end, label=right:
                    {
                    0.25~~~~~
                    $\begin{array}{l}
                    c_1:4\cdot 3+4\cdot 6\geq 30\\
                    c_2:4\cdot 3 = 12
                    \end{array}$
                    }] {}
                edge from parent
                node[above] {$s_2=4$}
                node[below]  {$0.5$}
            }
            child {
                node[end, label=right:
                    {
                    0.25~~~~~
                    $\begin{array}{l}
                    c_1:4\cdot 3+3\cdot 6\geq 30\\
                    c_2:3\cdot 3 \neq 12
                    \end{array}$
                    }] {}
                edge from parent
                node[above] {$s_2=3$}
                node[below]  {$0.5$}
            }
        edge from parent         
            node[above] {$s_1=4$}
            node[below]  {$0.5$}
    };
\end{tikzpicture}
\caption{Policy tree for the SCSP in Example 1}
\label{fig:msp_pt}
\end{figure}
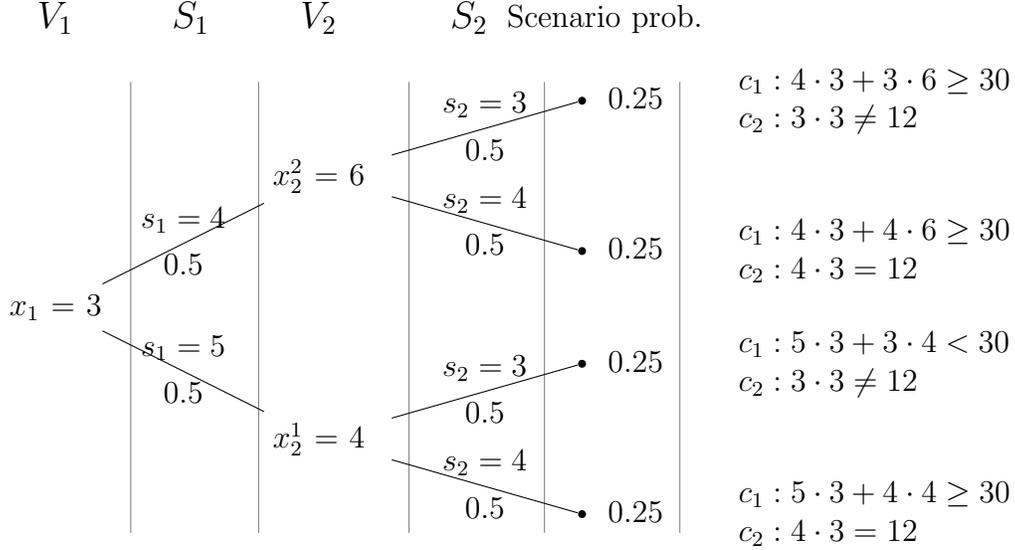

\begin{figure}[ht]
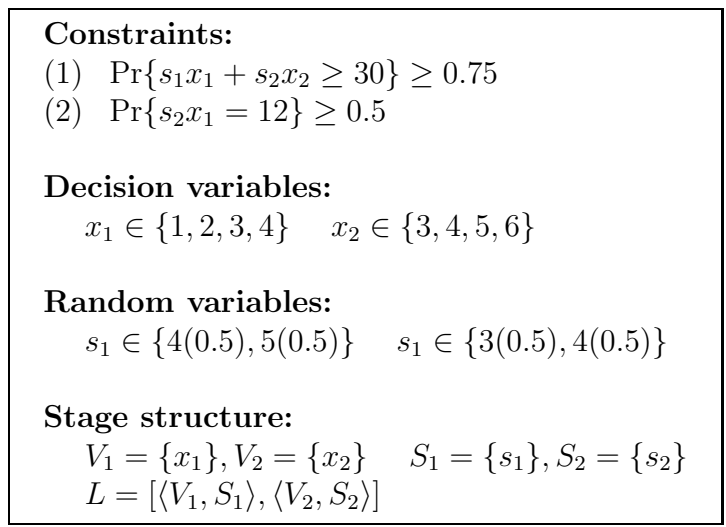

    \begin{center}
        \framebox{
        \begin{tabular}{ll}
        \mbox{{\bf Constraints:}} \\
        (1) $~\Pr\{s_1x_1+s_2x_2\geq 30\} \geq 0.75$\\
        (2) $~\Pr\{s_2x_1 = 12\} \geq 0.5$\\\\
        \mbox{{\bf Decision variables:}} \\
        $~~~~x_1 \in \{1,2,3,4\}~~~~x_2 \in \{3,4,5,6\}$ \\\\
        \mbox{{\bf Random variables:}} \\
        $~~~~s_1 \in \{4(0.5),5(0.5)\}~~~~s_1 \in \{3(0.5),4(0.5)\}$ \\\\
        \mbox{{\bf Stage structure:}} \\
        $~~~~V_1=\{x_1\}, V_2=\{x_2\}~~~~S_1=\{s_1\}, S_2=\{s_2\}$\\
        $~~~~L=[\langle V_1,S_1 \rangle,\langle V_2,S_2 \rangle]$
        \end{tabular}
        }
    \end{center}
    \caption{The two-stage SCSP in Example 1} \label{model:example_1_scsp}
\end{figure}

\begin{example}\label{example_1}
Let us consider the two-stage SCSP in Fig. \ref{model:example_1_scsp}, whose stage structure is
$L=[\langle V_1,S_1 \rangle,\langle V_2,S_2 \rangle]$;
$V_1=\{x_1\}$ and $S_1=\{s_1\}$, $V_2=\{x_2\}$ and $S_2=\{s_2\}$.
Random variable $s_1$ may take two possible values, 5 and 4, each with probability 0.5;
random variable $s_2$ may also take two possible values, 3 and 4, each with probability 0.5.
The domain of $x_1$ is $\{1,\ldots,4\}$, the domain of $x_2$ is $\{3,\ldots,6\}$. 
There are two chance constraints\footnote{In what follows, for convenience, we denote a 
chance constraint by using the notation ``$\Pr\{\langle cons\rangle\}\geq\beta$'', meaning that 
constraint $\langle cons\rangle$, constraining decision and random variables, should be satisfied
with probability greater than or equal to $\beta$.} 
in $C$,  $\Pr\{c_1: s_1x_1+s_2x_2\geq 30\} \geq
0.75 $ and $\Pr\{c_2: s_2x_1 = 12\} \geq
0.5 $. In this case, the decision variable $x_1$ must be set to a
unique value before random variables are observed, while decision
variable $x_2$ takes a value that depends on the observed value of
the random variable $s_1$. A possible solution to this SCSP is the
satisfying policy tree shown in Fig. \ref{fig:msp_pt} in which
$x_1=3,x_2^1=4$ and $x_2^2=6$, where $x_2^1$ is the value assigned
to decision variable $x_2$, if random variable $s_1$ takes value
$5$, and $x_2^2$ is the value assigned to decision variable $x_2$,
if random variable $s_1$ takes value $4$. 
\end{example}

\noindent
As Example \ref{example_1} shows, a solution to an SCSP is not simply an
assignment of the decision variables to values, but it is
instead a satisfying policy tree.

It is worth noting that {\em asking individual constraints to be satisfied according to their probability thresholds} is different from {\em asking a conjunction of constraints to be satisfied according to a prescribed probability threshold}. Informally speaking, in Example \ref{example_1} we simply state that $c_1: s_1x_1+s_2x_2\geq 30$ should hold true with probability $\beta_1=0.75$, i.e. in at least 75\% of the scenarios. If $c_2: s_2x_1 = 12$ holds true or not in those very same scenarios is not a matter of concern, as long as $c_2$ holds true in at least 50\% of the scenarios --- not necessarily the same as those in which $c_1$ holds true. Essentially, in Example \ref{example_1} we do not state anything about the conjunction $c_3: (s_1x_1+s_2x_2\geq 30) \wedge (s_2x_1 = 12)$. If we want to state something about this conjunction, we need to post a specific chance constraint $c_3$ with its own satisfaction threshold $\beta_3$. Assuming $\beta_3=0.5$, we may for instance require the conjunction $c_3$ to hold true in at least 50\% of the scenarios, i.e. $\Pr\{c_3: (s_1x_1+s_2x_2\geq 30) \wedge (s_2x_1 = 12)\}\geq 0.5$. Incidentally, the policy tree presented in Fig. \ref{fig:msp_pt} also satisfies this constraint. 

A practical example that further clarifies this discussion is found in inventory control. It is customary in inventory control problems to enforce service level constraints such as
\[\Pr\{I_t\geq 0\}\geq \alpha\hspace{3em}t=1,\ldots,N\]
where $N$ represents the length of the planning horizon and $I_t$ is the inventory level at the end of period $t$. The above set of constraints means that the probability of stocking out in a given period should be less than $1-\alpha$; regardless what happens in other periods. A more restrictive service level requirement would be 
\[\Pr\{\bigwedge_{t=1}^N I_t\geq 0\}\geq \alpha\]
This latter restriction means that the probability of stocking out in at least one of the $N$ periods should be less than $1-\alpha$.

\subsection{Confidence interval analysis}

Confidence interval analysis is a well established technique in statistics. 
Informally, confidence intervals are a useful tool for computing, from a 
given set of experimental results, a range of values 
that, with a certain confidence level (or confidence probability), 
will cover the actual value of a parameter that is being estimated.
Consider a discrete random variable that follows a Bernoulli distribution.
Accordingly, such a variable may produce only two outcomes, i.e. ``yes''
and ``no'', with probability $q$ and $1-q$, respectively. 
Let us assume that the value $q$ --- the ``yes'' probability --- is unknown. 
Obviously, if we observe the outcome of a Bernoulli trial once, the data collected will not reveal 
much about the value of $q$. Nevertheless,
in practice, we may be interested in ``estimating'' $q$, by repeatedly 
observing the behavior of the random variable in a sequence of
Bernoulli trials. 
This problem is well-known in 
statistics and both exact and approximate techniques are available
for performing this estimation \cite{cloppears34,agresticloull98}.
The estimation produced by the methods available in the literature typically
does not come as a point estimate, rather it consists of an interval of values
computed from a set of representative samples for the quantity being estimated.
This interval is known as ``confidence interval'' and consists of a range of values 
that, with a certain confidence probability $\alpha$, 
covers the actual value of the parameter that is being estimated.

A method that is commonly classified as the ``exact confidence intervals'' for the 
Binomial distribution has been introduced by Clopper and Pearson in \cite{cloppears34}.
This method uses the Binomial cumulative distribution function (CDF) in order to build the interval
from the data observed. The Clopper-Pearson interval is a symmetric two-sided confidence interval. It can be however also expressed as a single-sided interval. Clopper-Pearson single-sided intervals can be written as
$\left(p_{\text{lb}},1\right)$ and $\left(0,p_{\text{ub}}\right)$ where\\
\[
\begin{array}{ll}
&p_{\text{lb}}=\min\{q| \Pr\{\text{bin}(N;q)\geq X\}\geq 1-\alpha\},\\
&p_{\text{ub}}=\max\{q| \Pr\{\text{bin}(N;q)\leq X\}\geq 1-\alpha\},
\end{array}
\]
$X$ is the number of successes (or ``yes'' events) observed in the sample, $\text{bin}(N; q)$ is a binomial random variable with $N$ trials and probability of success $q$ and $\alpha$ is the confidence probability. Note that we assume $p_{\text{lb}}=0$ when $X=0$ and that $p_{\text{ub}}=1$ when $X=N$.

Because of the close relationship between Binomial distribution and the Beta distribution, 
the Clopper-Pearson interval is sometimes presented in an alternative format 
that uses percentiles from the beta distribution \cite{ehp00}:
\[
\begin{array}{ll}
&p_{\text{lb}}=1-\text{beta}^{-1}(\alpha,N-X+1,X),\\
&p_{\text{ub}}=1-\text{beta}^{-1}(1-\alpha,N-X,X+1),
\end{array}
\]
where $\text{beta}^{-1}$ denotes the inverse Beta distribution.
This form can be efficiently evaluated by existing algorithms.

An interesting property of confidence intervals related to the estimation of
the ``success'' probability associated with a Bernoulli trial
consists in the fact that, given a confidence probability, it is possible to derive mathematically,
by performing a worst case analysis, the minimum number of samples that should be
observed in order to produce a confidence interval of a given size.

Therefore, for a given confidence probability $\alpha$, it is possible to determine
the minimum number of samples
that should be considered in order to achieve a margin of error of $\pm\vartheta$ 
in the estimation of the ``success'' probability of a Bernoulli trial. 
This computation plays a central role in our novel approach.
In fact, intuitively estimating the satisfaction probability of a chance constraint
is equivalent to estimating the ``success'' probability of the associated Bernoulli trial.

\section{Sampled SCSPs}\label{sec:sampled_scsps}

Consider an SCSP $\mathcal{P}$ over a set $S$ of stochastic variables. Assume that
stochastic variables are defined on continuous supports or 
discrete supports comprising a large or infinite number of values. 
Solving the original SCSP clearly poses a hard combinatorial challenge, in fact
the policy tree comprises a number of scenarios that is exponential in the size of
stochastic variable domains. Since the policy tree may comprise an infinite number of scenarios, 
the computational problem may even be undecidable in general.

In this section we discuss how to {\em sample} a more compact SCSP, which 
comprises at most $N$ scenarios, out of the original problem. 
We shall call this new problem $\widehat{\mathcal{P}}_N$ or ``sampled SCSP'' over $N$
scenarios. Intuitively, a sampled SCSP is a reduced version of the original problem 
the solution of which is a policy tree that comprises a bounded number of 
paths sampled out of the original policy tree.
In the following sections we will discuss under which conditions the solution
to a sampled SCSP $\widehat{\mathcal{P}}_N$ is, according to a certain confidence 
probability and within prescribed error tolerance thresholds, likely to be also a solution to the original SCSP $\mathcal{P}$.

We shall here discuss how to employ Simple Random Sampling \cite{uc02} to obtain
a sampled SCSP from the original problem. Of course, more advanced stratified
sampling techniques may be used in order to reduce variance and improve the 
effectiveness of the approach. Nevertheless, due to the large number of topics
already covered in this work, we leave this discussion as future work.

Consider a complete realization, $\bar{s}_1,\ldots,\bar{s}_{m}$, for the 
stochastic variables in $S$ obtained by sampling a value from the support $D(s_i)$ of each of the stochastic
variables $s_i\in S$ according to its probability distribution $P(s_i)$. From the definition of policy tree
it is clear that there always exists a path associated with
this realization. In other words, this realization corresponds to one of the scenarios
comprised in the policy tree.

Generate $N$ independent sets of random variable realizations 
\[\{\bar{s}^1_1,\ldots,\bar{s}^1_{m}\},\{
\bar{s}^2_1,\ldots,\bar{s}^2_{m}\},\ldots,
\{\bar{s}^N_1,\ldots,\bar{s}^N_{m}\},\]
where $\bar{s}^i_j$ is the realized value for random variable $j$ observed in the $i$-th 
set of realizations. 
Recall that $\mathcal{T}$ denotes the complete, and possibly infinite, policy tree for $\mathcal{P}$.
Let a reduced policy tree $\widehat{\mathcal{T}}$ for $\mathcal{P}$ be a policy tree that comprises only 
arc labelings observed in the former $N$ complete realizations (without repetitions).

Let $\widehat{\Psi}$ denote the reduced set of distinct paths in  $\widehat{\mathcal{T}}$.
The probability associated with each path $p\in \widehat{\Psi}$, i.e. 
$\Pr\{\text{arcs}(p)\}$, is simply set equal to the frequency 
of occurrence of such a path in the above $N$ realizations.
Of course, $\widehat{\mathcal{T}}$ represents a policy tree for a different SCSP
than the one we started with. We call this new problem the 
sampled SCSP $\widehat{\mathcal{P}}_N$.

Now consider a chance constraint $h \in C$ with a specified
threshold level $\beta_h$, a policy tree $\widehat{\mathcal{T}}$
for the sampled SCSP $\widehat{\mathcal{P}}_N$ and a path $p \in \mathcal{T}$. 
We say that $h$ is {\em satisfied wrt to a given policy tree $\widehat{\mathcal{T}}$} iff
\[\sum_{p \in \widehat{\Psi}: \bar{h}_{\downarrow{p}} \in h_{\downarrow{p}}} \Pr\{\text{arcs}(p)\} \geq
\beta_h.\]
\begin{example}\label{example_2}
Let us consider the two-stage SCSP $\mathcal{P}$ discussed in Example \ref{example_1}. We set
$N=3$ and we derive a sampled SCSP $\widehat{\mathcal{P}}_N$. By using simple
random sampling we draw the following three complete realizations for random 
variables in $\mathcal{P}$:
\[\{\bar{s}^1_1=5,\bar{s}^1_2=4\},
\{\bar{s}^2_1=4,\bar{s}^2_2=4\},
\{\bar{s}^3_1=5,\bar{s}^3_2=4\}.\]
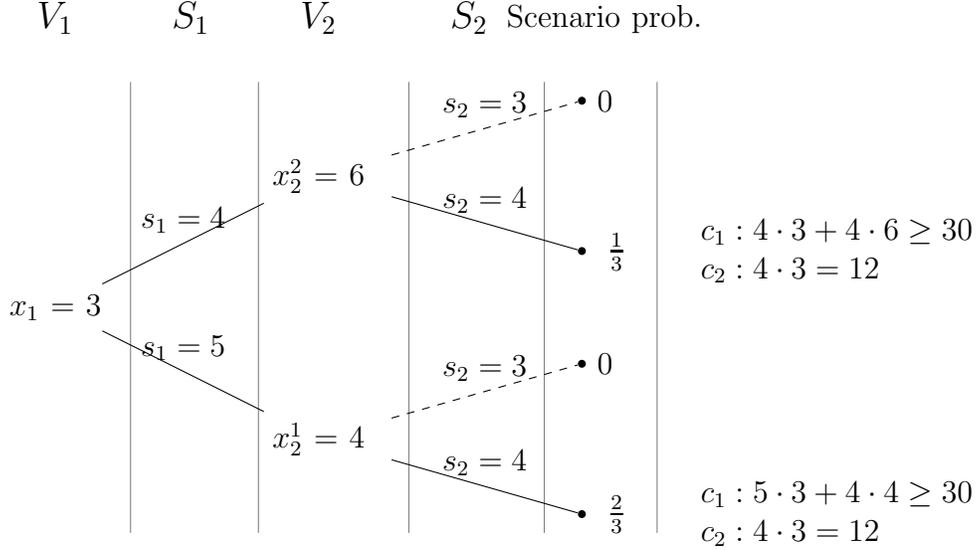
\begin{figure}[t!]
\centering 
\tikzstyle{level 1}=[level distance=3.5cm, sibling distance=3.5cm]
\tikzstyle{level 2}=[level distance=3.5cm, sibling distance=2cm]

\tikzstyle{bag} = [text width=4em, text centered]
\tikzstyle{end} = [circle, minimum width=3pt,fill, inner sep=0pt]


\begin{tikzpicture}[grow=right]
\draw[name path=V1a--V1b,style=thin,color=gray] (1,-3) -- (1,3);
\draw[name path=V1a--V1b,style=thin,color=gray] (2.7,-3) -- (2.7,3);
\draw[name path=V1a--V1b,style=thin,color=gray] (4.7,-3) -- (4.7,3);
\draw[name path=V1a--V1b,style=thin,color=gray] (6.5,-3) -- (6.5,3);
\draw[name path=V1a--V1b,style=thin,color=gray] (8,-3) -- (8,3);
\draw (0,3.5) node[above] {\large$V_1$};
\draw (1.8,3.5) node[above] {\large$S_1$};
\draw (3.5,3.5) node[above] {\large$V_2$};
\draw (5.5,3.5) node[above] {\large$S_2$};
\draw (7.3,3.5) node[above] {Scenario prob.};
\node[bag] {$x_1=3$}
    child {
        node[bag] {$x_2^1=4$}        
            child {
                node[end, label=right:
                    {
                    $\frac{2}{3}$~~~~~
                    $\begin{array}{l}
                    c_1:5\cdot 3+4\cdot 4\geq 30\\
                    c_2:4\cdot 3 = 12
                    \end{array}$
                    }] {}
                edge from parent
                node[above] {$s_2=4$}
            }
            child {
                node[end, label=right:
                    {0}] {}
                edge from parent[style=dashed]
                node[above] {$s_2=3$}
            }
            edge from parent 
            node[above] {$s_1=5$}
    }
    child {
        node[bag] {$x_2^2=6$}        
        child {
                node[end, label=right:
                    {
                    $\frac{1}{3}$~~~~~
                    $\begin{array}{l}
                    c_1:4\cdot 3+4\cdot 6\geq 30\\
                    c_2:4\cdot 3 = 12
                    \end{array}$
                    }] {}
                edge from parent
                node[above] {$s_2=4$}
            }
            child {
                node[end, label=right:
                    {0}] {}
                edge from parent[style=dashed]
                node[above] {$s_2=3$}
            }
        edge from parent         
            node[above] {$s_1=4$}
    };
\end{tikzpicture}
\caption{Reduced policy tree for the sampled SCSP in Example \ref{example_2}; dashed arcs are those that have not been observed in the sample.}
\label{fig:msp_pt_sampling}
\end{figure}
A possible solution to the sampled SCSP $\widehat{\mathcal{P}}_N$ is the
satisfying policy tree shown in Fig. \ref{fig:msp_pt_sampling}, in which
$x_1=3,x_2^1=4$ and $x_2^2=6$, where $x_2^1$ is the value assigned
to decision variable $x_2$, if stochastic variable $s_1$ takes value
$5$, and $x_2^2$ is the value assigned to decision variable $x_2$,
if stochastic variable $s_1$ takes value $4$. The 
above policy tree has two paths sampled out of the original tree:
$p_1$ has an associated probability of 2/3, since we observed two occurrences
of the scenario associated with this path over the 3 complete realisations 
sampled for the random variables;
$p_2$ has an associated probability of 1/3, since we observed a single occurrence
of the scenario associated with this path over the 3 complete realisations sampled for 
the random variables.
Paths that were not observed in the sampled realisations have an associated probability 
equal to zero and are not considered. 
\end{example}
It should be noted that every policy tree $\widehat{\mathcal{T}}$ for a sampled SCSP
$\widehat{\mathcal{P}}$ can be employed as a (partial) policy tree for the original SCSP
$\mathcal{P}$. 
Nevertheless, by sampling we lose completeness. 
If at stage $i$ in $\mathcal{P}$ we observe, for a given random variable, a realised value  
that is not comprised in $\widehat{\mathcal{T}}$, it will be of course impossible to determine
the correct decisions for subsequent stages. By taking a conservative point of view, this means that all paths in the corresponding subtree will never be satisfied. 
In multi-stage SCSPs, and especially in those including random variables with continuous support,
this prevents the direct use of the approach that will be discussed in this work. In fact, if random variable supports are continuous, there is only an infinitesimal probability of observing a given set of realisations.
In this case, it is therefore essential to adopt a ``rolling horizon'' approach \cite{115085} in order
to reduce the original multi-stage SCSPs to a sequence of multi-stage sampled SCSPs. 
Under this strategy, our aim is to fix decisions at stage one, and make sure that compatible
values exist for decision variables that appear, for subsequent stages, in $\widehat{\mathcal{T}}$.
Future decisions are not fixed because, after observing the realised values for 
random variables at stage one, the problem is solved again by taking into account
new available information; decision variables that were previously 
associated with stage two become stage one decisions.
The original problem is thus reduced to a sequence of multi-stage sampled SCSPs.
We will apply this technique to handle the two-stage problem 
discussed in Section \ref{sec:smps}: the stochastic multiprocessor 
scheduling problem with release time and deadlines.

\section{($\alpha$,$\vartheta$)-solutions}\label{sec:sampled_scsps_solutions}

We will now characterize the probability that the solution of a sampled SCSPs 
$\widehat{\mathcal{P}}_N$ over $N$ scenarios, which may be computed by using any of
the existing approaches discussed in \cite{tmw2006,citeulike:10689458}, is a solution 
to the original {\em single-stage} SCSP $\mathcal{P}$. 

These results
are also applicable to multi-stage problems, provided that a rolling horizon approach is adopted and
that the aim is to characterize the probability that stage one decisions of a sampled SCSPs 
$\widehat{\mathcal{P}}_N$ over $N$ scenarios are consistent with respect to
the original SCSP $\mathcal{P}$.

We will firstly discuss how to compute $N$ such that, if a given policy tree $\mathcal{T}$ 
satisfies a chance constraint $h$ in the sampled SCSPs $\widehat{\mathcal{P}}_N$, it also 
satisfies the same chance constraint in the original SCSP $\mathcal{P}$ with probability $\alpha$.
Since a policy tree $\mathcal{T}$ in $\widehat{\mathcal{P}}_N$
by definition only comprises a subset $\widehat{\Psi}$ of all the paths that constitute a policy tree
for the original SCSP $\mathcal{P}$, this policy tree, in order to satisfy $h$ in the original SCSP 
$\mathcal{P}$, must clearly provide a sufficient satisfaction probability regardless of the scenarios
that have been ignored by the sampling process.

Consider a confidence probability $\alpha$ and a margin of error of $\pm\vartheta$;
The number of scenarios $N$ for the sampled SCSP depends on $\vartheta$, $\alpha$ and also 
$\beta$, which we recall is the target satisfaction probability of chance constraint $h$.

\begin{definition}\label{sample_size}
$N$ is computed as the minimum value for which \[\max(p^{\beta}_{\text{ub}}-\beta,\beta-p^{\beta}_{\text{lb}})\leq \vartheta,\]
where $p^{\beta}_{\text{lb}}$ and $p^{\beta}_{\text{ub}}$ are the single-sided Clopper-Pearson confidence interval bounds for a confidence probability $\alpha$, and $\mbox{\it round}(\beta N)$ ``successes'' in $N$ trials; $\mbox{\it round}()$ approximates the value to the nearest integer.\footnote{This is justified by the fact that the Clopper-Pearson interval is, in fact, a step function --- see \cite{cloppears34}, p. 405 --- since the Binomial is a discrete probability distribution.}
\end{definition}

\begin{definition} \label{consistency_condition_beta}
Any policy tree $\mathcal{T}$, which can be proved to satisfy $h$ in $\mathcal{P}$ with probability $\alpha$,
satisfies $h$ in $\mathcal{P}$ with probability $\alpha$ if it satisfies $h$ in $\widehat{\mathcal{P}}_N$. Conversely,
any policy tree $\mathcal{T}$, which can be proved not to satisfy $h$ in $\mathcal{P}$ with probability $\alpha$,
does not satisfy $h$ in $\mathcal{P}$ with probability $\alpha$, if it does not satisfy $h$ in $\widehat{\mathcal{P}}_N$.
\end{definition}

\begin{prop} \label{proposition1}
A policy tree $\mathcal{T}$ can be proved to satisfy $h$ in $\mathcal{P}$ with probability $\alpha$ if 
the actual satisfaction probability $\delta>\beta$ provided by $\mathcal{T}$ wrt $h$ is such that $\delta\geq p^{\beta}_{\text{ub}}$. Conversely, if the actual satisfaction probability $\delta<\beta$ provided by $\mathcal{T}$ wrt $h$ is such that $\delta\leq p^{\beta}_{\text{lb}}$ $\mathcal{T}$ can be proved to not satisfy $h$ in $\mathcal{P}$ with probability $\alpha$.
\end{prop}
\begin{proof}
Let $\delta\geq p^{\beta}_{\text{ub}}$. Since 
$p^{\beta}_{\text{ub}}=\max\{q| \Pr\{\text{bin}(N;q)\leq \text{round}(\beta N)\}\geq 1-\alpha$,
it is clear that $\Pr\{\text{bin}(N;\delta)\leq \text{round}(\beta N)\}< 1-\alpha$. This
means that \[\Pr\left\{\sum_{p \in \widehat{\Psi}: \bar{h}_{\downarrow{p}} \in h_{\downarrow{p}}} \Pr\{\text{arcs}(p)\}\leq\beta\right\}<1-\alpha,\] where we recall that $\widehat{\Psi}$ is the set of paths in the sampled SCSP $\widehat{\mathcal{P}}_N$.
This implies \[\Pr\left\{\sum_{p \in \widehat{\Psi}: \bar{h}_{\downarrow{p}} \in h_{\downarrow{p}}} \Pr\{\text{arcs}(p)\}\geq\beta\right\}\geq\alpha.\] Therefore, by using the test \[\sum_{p \in \widehat{\Psi}: \bar{h}_{\downarrow{p}} \in h_{\downarrow{p}}} \Pr\{\text{arcs}(p)\}\geq\beta,\] a policy tree $\mathcal{T}$ can be proved to satisfy $h$ in $\mathcal{P}$ with probability $\alpha$. 

Conversely, let $\delta\leq p^{\beta}_{\text{lb}}$. Since 
$p^{\beta}_{\text{lb}}=\min\{q| \Pr\{\text{bin}(N;q)\geq \text{round}(\beta N)\}\geq  1-\alpha$,
it is clear that \[\Pr\{\text{bin}(N;\delta)\geq \text{round}(\beta N)\}< 1-\alpha.\] This
means that \[\Pr\left\{\sum_{p \in \widehat{\Psi}: \bar{h}_{\downarrow{p}} \in h_{\downarrow{p}}} \Pr\{\text{arcs}(p)\}\geq\beta\right\}<1-\alpha,\] which implies \[\Pr\left\{\sum_{p \in \widehat{\Psi}: \bar{h}_{\downarrow{p}} \in h_{\downarrow{p}}} \Pr\{\text{arcs}(p)\}\leq\beta\right\}\geq\alpha.\] Therefore, by using the test \[\sum_{p \in \widehat{\Psi}: \bar{h}_{\downarrow{p}} \in h_{\downarrow{p}}} \Pr\{\text{arcs}(p)\}\leq\beta,\] a policy tree $\mathcal{T}$ can be proved to not satisfy $h$ in $\mathcal{P}$ with probability $\alpha$. 
\end{proof}

\begin{prop} \label{proposition2}
Any policy tree $\mathcal{T}$ which provides a satisfaction probability $\delta\geq \beta+\vartheta$
wrt $h$ in $\mathcal{P}$ can be proved to satisfy $h$ in $\mathcal{P}$ with probability $\alpha$.
Any policy tree $\mathcal{T}$ which provides a satisfaction probability $\delta\leq \beta-\vartheta$
wrt $h$ in $\mathcal{P}$ can be proved to not satisfy $h$ in $\mathcal{P}$ with probability $\alpha$.
\end{prop}
\begin{proof}
this directly follows from Definition \ref{sample_size} and Proposition \ref{proposition1}.
\end{proof}

\begin{prop} \label{proposition3}
Any policy tree $\mathcal{T}$ which can not be proved to satisfy or not to satisfy $h$ in $\mathcal{P}$ with probability $\alpha$, can be either proved to satisfy $h$ in $\mathcal{P}$ with probability $\gamma$, where $\gamma$ is a probability ranging in $[0.5,\alpha)$, if it satisfies $h$ in $\widehat{\mathcal{P}}_N$, or not to satisfy $h$ in $\mathcal{P}$ with probability $\gamma$, where $\gamma$ is a probability ranging in $[0.5,\alpha)$, if it does not satisfies $h$ in $\widehat{\mathcal{P}}_N$.
\end{prop}
\begin{proof}
Consider the two limiting cases. (i) The actual satisfaction probability $\delta$ provided by $\mathcal{T}$ wrt $h$ in $\mathcal{P}$ is exactly equal to $\beta$. Since the sample mean, used to estimate the satisfaction probability out of the $N$ samples considered, is an unbiased estimator of $\delta$, it will overestimate $\beta$ with probability $0.5$ and, similarly, it will underestimate $\beta$ with probability $0.5$; this sets the lower bound for $\gamma$. (ii) The actual satisfaction probability $\delta$ provided by $\mathcal{T}$ wrt $h$ in $\mathcal{P}$ is exactly equal to $\beta+\vartheta$. From the proof of Proposition \ref{proposition1} it immediately follows that, in this case, $\gamma=\alpha$, and also that, if $\delta<\beta+\vartheta$ then $\gamma<\alpha$; this sets the upper bound for $\gamma$.
\end{proof}

\begin{definition}\label{alpha_theta_solution}
An $(\alpha,\vartheta)$-solution to an SCSP $\mathcal{P}$ is a policy tree $\widehat{\mathcal{T}}$
that at least with probability $\alpha$ provides
for every chance constraint $h_i$ in $\mathcal{P}$ with satisfaction threshold $\beta_i$ a 
satisfaction probability greater than or equal to $\beta_i-\vartheta$.
\end{definition}

It is apparent that $\vartheta$ may be interpreted as a parameter that the user can set in order to define a ``region of indifference'', i.e. $\beta\pm \vartheta$, for the satisfaction probability. In such a region, we assume that assignments can be safely misclassified with probability greater than $\alpha$ and that satisfaction probabilities remain in an acceptable range.

\begin{example}\label{example3}
Consider the single-stage SCSP $\mathcal{P}=\langle
V,S,D,P,C,\beta, L \rangle$, where
$V=\{X_1,X_2\}$, $S=\{r_1,r_2\}$, $D(X_1)=D(X_2)=\{0,1\}$,  
$D(r_1)=(0,100)$, $P(r_1)=\text{\normalfont uniform}(0,100)$, 
$D(r_2)=(0,300)$, $P(r_2)=\text{\normalfont uniform}(0,300)$, 
$C=\{c:C_1\geq X_1 r_1 + X_2 r_2\}$, $\beta_c=0.5$, and $L=[\langle V,S\rangle]$.
$C_1=185$ is a constant. This problem comprises random variables
defined on a continuous support and it cannot be solved by existing complete approaches
to SCSPs. If we set $\alpha=0.95$ and $\vartheta=0.05$, from Definition \ref{sample_size} we compute
the number of samples $N=290$ required to guarantee that any solution to the sampled SCSP
$\widehat{\mathcal{P}}$ over $N$ samples is an $(\alpha,\vartheta)$-solution for $\mathcal{P}$. 

Furthermore, the simple structure of the constraint $c$ considered in $\mathcal{P}$ allows us to perform some further analysis. Consider the assignment $X_1=1$ and $X_2=1$. A simple reasoning on the convolution of two independently non-identically distributed uniform random variables (see \cite{RePEc:spr:stpapr:v:50:y:2009:i:1:p:171-175}) immediately suggests that this assignment is indeed inconsistent.
\begin{figure}[ht]
\centering
\begin{minipage}[c]{0.475\textwidth}
\centering
\includegraphics[type=eps,ext=.eps,read=.eps,width=1\columnwidth]{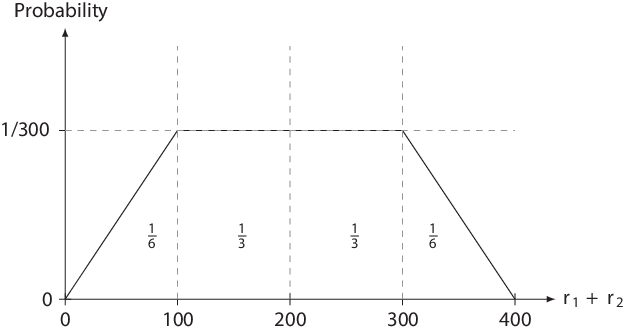}
\caption{Probability density function of the convolution of two independently non-identically distributed uniform random variables $r_1$ and $r_2$.}
\label{fig:convolution_uniform}
\end{minipage}
\hspace{0.025\linewidth}
\begin{minipage}[c]{0.475\textwidth}
\includegraphics[type=eps,ext=.eps,read=.eps,width=1\columnwidth]{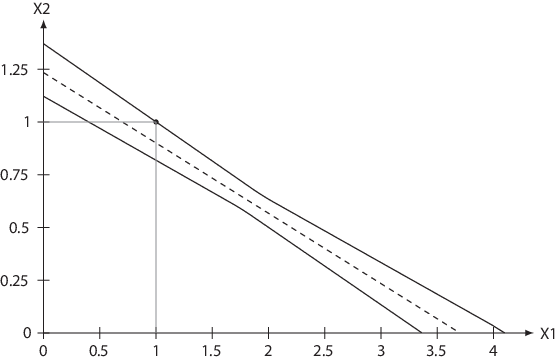}
\caption{Feasible region for the SCSP in Example 1; the dashed line denotes the true boundary of  constraint $c$. The upper solid line demarcates the set of solutions providing a satisfaction probability of at least $\beta-\vartheta$; the lower solid line demarcate the set of solutions providing a satisfaction probability of at least $\beta+\vartheta$.}
\label{fig:feasible_region_example1}
\end{minipage}
\end{figure}
$r_1$ and $r_2$ are two independently non-identically distributed uniform random variables. 
The distribution that results from their convolution is shown in Fig. \ref{fig:convolution_uniform}. This distribution is shaped like a trapezoid. Clearly, since the area for the whole figure must be equal to 1, the area of each of the two rectangle triangles at the side of the trapezoid must be equal to 1/6. Consequently, the area of the internal rectangle must be equal to 2/3. It is easy to see that the cumulative distribution function for value 200 returns a probability of 0.5. Then, since 1/3*(15/100)=0.05, the 0.45 quantile of the  inverse cumulative distribution function which results from convoluting $r_1$ and $r_2$ is exactly equal to $C_1=185$. Therefore, since the satisfaction probability provided by the assignment $X_1=1$ and $X_2=1$ is equal to $\beta_c-\vartheta=0.45$ (Fig. \ref{fig:feasible_region_example1}), this assignment will be correctly classified as inconsistent with probability $\alpha$, when the sample size is set to $N=290$. 
\end{example}

Let $h_1,\ldots, h_k$ be $k$ chance constraints in an SCSP $\mathcal{P}$.
Let  $\widehat{\mathcal{P}}$ be a sampled SCSP over $N$ samples,
where $N$ is the number of samples required to guarantee a confidence 
level $\alpha$ and an error tolerance threshold $\vartheta$ for each 
constraint $h_i$ considered independently, according to Definition \ref{sample_size}.
\begin{prop} \label{proposition4}
Let $\widehat{\mathcal{T}}$ be a policy tree that is a solution to 
$\widehat{\mathcal{P}}$. Then $\widehat{\mathcal{T}}$
is an $(\alpha,\vartheta)$-solution for $\mathcal{P}$.
\end{prop}
\begin{proof}
Consider a chance constraint $h_i$. Let $\beta_i$ be the respective
satisfaction threshold. By definition, the probability that a solution $\widehat{\mathcal{T}}$
to $\widehat{\mathcal{P}}$ provides a service level less or
equal to $\beta_i-\vartheta$ for $h_i$ in $\mathcal{P}$ is less than or equal to $1-\alpha$.
Therefore $\widehat{\mathcal{T}}$ is an $(\alpha,\vartheta)$-solution.
Now consider a pair of chance constraints $\langle h_i, h_j\rangle$ with 
satisfaction thresholds  $\beta_i$,  $\beta_j$, respectively.
The probability $p_{ij}$ that a solution $\widehat{\mathcal{T}}$
to $\widehat{\mathcal{P}}$ provides a service level less or
equal to $\beta_i-\vartheta$ for $h_i$ and to $\beta_j-\vartheta$ for 
$h_j$ in $\mathcal{P}$ is less than or equal to $(1-\alpha)^2$, in fact
we must misclassify both the constraints in order to accept such a solution.
Even a single constraint correctly classified will make $\widehat{\mathcal{T}}$
inconsistent w.r.t. $\widehat{\mathcal{P}}$. 
The case in which constraints are misclassified independently from each other represents 
a worst-case reasoning. If constraints are perfectly positively
correlated, i.e. if one is misclassified then all other constraints are also misclassified, then 
$p_{ij}$ is $(1-\alpha)$; if constraints are perfectly negatively correlated, i.e. if one is 
misclassified then no other constraint is misclassified, $p_{ij}$ becomes 0.
This reasoning can be generalized to $k$ chance constraints, 
for which the probability becomes $(1-\alpha)^k$. 
Noting that $(1-\alpha)^k<\ldots<(1-\alpha)^2<(1-\alpha)$ and that $1-(1-\alpha)^k\geq \alpha$ 
the probability $1-\alpha$ that a solution is misclassified in a model comprising
a single constraint represents an upper bound for the 
probability that a solution $\widehat{\mathcal{T}}$ to $\widehat{\mathcal{P}}$
does not provide a satisfaction probability within the required
tolerance threshold for one or more constraints in a generic model $\mathcal{P}$. By rephrasing,
the probability that a solution $\widehat{\mathcal{T}}$
provides a satisfaction probability greater than or equal to 
$\beta_i-\vartheta$ for each constraint $h_i$ is greater than or equal to $\alpha$. 
Hence, by Definition \ref{alpha_theta_solution}, $\widehat{\mathcal{T}}$ 
is an $(\alpha,\vartheta)$-solution for $\mathcal{P}$.  
\end{proof}

\section{($\alpha$,$\vartheta$)-solution set}\label{sec:alpha_theta_solution_set}

Consider policy tree $\mathcal{T}$, chance constraint $h$, and the indicator random variable
\[\tau= \left\{
\begin{array}{ll}
1&\sum_{p \in \widehat{\Psi}: \bar{h}_{\downarrow{p}} \in h_{\downarrow{p}}} \Pr\{\text{arcs}(p)\}\geq\beta\\
0&\text{otherwise}
\end{array}\right.
\] 
representing the test discussed in Definition \ref{consistency_condition_beta}. If the actual satisfaction probability $\delta$ provided by a policy tree $\mathcal{T}$ with respect to constraint $h$ in $\mathcal{P}$ is exactly equal to $\beta-\vartheta$, $\tau$ takes value 1 with probability $1-\alpha$.

To motivate the following discussion, we introduce the following example. 
\begin{example}
Consider once more Example \ref{example3} and 
assume that $D(X_1)=D(X_2)=(0,5)$; i.e. decision variables are defined on continuous domains spanning from 0 to 5. Assignments $(X_1=4.1,X_2=0)$ and $(X_1=0,X_2=1.37)$ lie on the upper solid line shown in Fig. \ref{fig:feasible_region_example1}. Each of these two assignments provides a satisfaction probability of exactly $\beta-\vartheta$ with respect to constraint $c$ in the original problem $\mathcal{P}$. From the discussion in Section \ref{sec:sampled_scsps_solutions} it follows that each of these two assignments is recognised as infeasible with probability $\alpha$ if $N=290$. However, since $r_1$ and $r_1$ are independent the probability that these two assignments are {\bf both} recognised as infeasible is only $\alpha^2$. We next discuss how to address the issue of correctly classifying multiple policy trees according to a prescribed confidence level $\alpha$.
\end{example}

We introduce the following definition.
\begin{definition}\label{alpha_theta_solution_set}
An $(\alpha,\vartheta)$-solution set to an SCSP $\mathcal{P}$ is a set of policy trees. All policy trees in this set simultaneously provide, with probability at least $\alpha$, a satisfaction probability greater than or equal to $\beta_i-\vartheta$ for every chance constraint $h_i$ in $\mathcal{P}$ with satisfaction threshold $\beta_i$.
\end{definition}

Consider an SCSP and $T$ policy trees $\mathcal{T}_1,\ldots,\mathcal{T}_T$ for which the actual satisfaction probability $\delta$ with respect to $h$ in $\mathcal{P}$ is less than or equal to $\beta-\vartheta$. Let $\tau_1\ldots,\tau_T$ be the associated random variables, each of which according to Proposition \ref{proposition4} takes value 1 with probability less than or equal to $1-\alpha$. Although we have fully characterised the marginal probability distribution of a test $\tau_i$ involving a single policy tree $\mathcal{T}_i$, we have not characterised yet the joint probability among tests carried out on a set of $T$ policy trees. 

\begin{prop}\label{prop:failure_probability_T}
The probability that $\tau_1,\ldots,\tau_T$ are all equal to 0 is at least $1 - T(1-\alpha)$.
\end{prop}
\begin{proof}
A worst-case reasoning can be carried out by considering the case in which events $\tau_i=1$ and $\tau_j=1$ are mutually exclusive for all $i,j=1,\ldots,T$, $i\neq j$; of course it is still true that $\Pr\{\tau_i=1\}=\Pr\{\tau_j=1\}\leq1-\alpha$. 
The probability that $\tau_1,\ldots,\tau_T$ are all equal to 0 is then easily seen to be $1-T(1-\alpha)$. If events are not mutually exclusive, this probability is greater than or equal to $1-T(1-\alpha)$, e.g. in the case of $T$ independent tests it would be $1-(1-\alpha)^T\geq 1-T(1-\alpha)$.
\end{proof}

Of course, it is not possible to know the value of $T$ a priori, as this would require solving the SCSP. However, for a given chance constraint $h$, $T$ is clearly less than or equal to the cardinality $A_h$ of the assignment space constrained by $h$. $A_h$ can be computed as the cartesian product of the domains of the decision variables in the policy tree that are constrained by $h$. Since the property discussed in Proposition \ref{prop:failure_probability_T} applies to each chance constraint $h\in C$, to compute an $(\alpha,\vartheta)$-solution set we may introduce the following Bonferroni's correction \cite{citeulike:13411582}, which is {\em free of correlation and distribution assumptions}, while computing $N$.

\begin{definition}\label{sample_size_alpha_theta_solution_set}
$N$ is computed as the minimum value for which \[\max(p^{\beta}_{\text{ub}}-\beta,\beta-p^{\beta}_{\text{lb}})\leq \vartheta,\]
where $p^{\beta}_{\text{lb}}$ and $p^{\beta}_{\text{ub}}$ are the single-sided Clopper-Pearson confidence interval bounds for a confidence probability $\widehat{\alpha}$, where  \[\widehat{\alpha}= 1-\frac{1-\alpha}{\sum_{h\in C}A_h},\] and
$\mbox{\it round}(\beta N)$ ``successes'' in $N$ trials.
\end{definition}

\begin{prop} \label{prop:bonferroni}
A set of policy trees that are solutions to $\widehat{\mathcal{P}}$ for a sample size $N$ computed as discussed in Definition \ref{sample_size_alpha_theta_solution_set} is an $(\alpha,\vartheta)$-solution set for $\mathcal{P}$.
\end{prop}
\begin{proof}
Bonferroni's correction, introduced in Definition \ref{sample_size_alpha_theta_solution_set}, ensures that, for every chance constraint $h$ in $C$, the probability $\tau_1,\ldots,\tau_T$ are all equal to 0 simultaneously is at least $\alpha$.
\end{proof}

\begin{example}\label{ex:example_5}
\begin{figure}[ht]
\centering
\begin{minipage}[t]{0.475\textwidth}
\centering
\includegraphics[type=eps,ext=.eps,read=.eps,width=0.98\columnwidth]{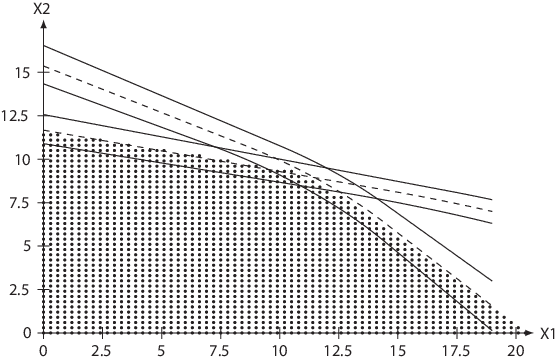}
\caption{An $(\alpha,\vartheta)$-solution set for Example \ref{ex:example_5} computed for $N=2848$ samples}
\label{fig:graph_region_2848}
\end{minipage}
\hspace{0.025\linewidth}
\begin{minipage}[t]{0.475\textwidth}
\includegraphics[type=eps,ext=.eps,read=.eps,width=1\columnwidth]{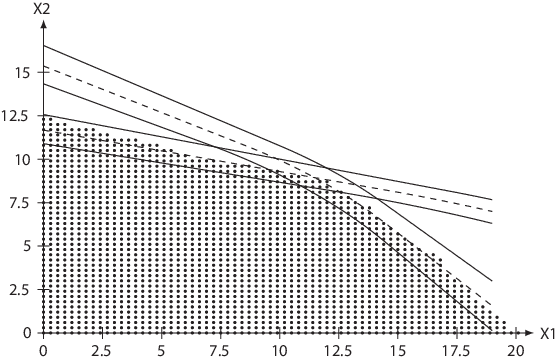}
\caption{An approximate $(\alpha,\vartheta)$-solution set for Example \ref{ex:example_6} computed for $N=348$ samples}
\label{fig:graph_region_348}
\end{minipage}
\end{figure}
Consider the following SCSP $\mathcal{P}=\langle
V,S,D,P,C,\beta_c, L \rangle$, where
$V=\{X_1,X_2\}$, $S=\{r_1,r_2\}$, $D(X_1)=D(X_2)=\{0,0.01,0.02,\ldots,24.99,25\}$,  
$D(r_1)=(0,10)$, $P(r_1)=\text{\normalfont uniform}(0,10)$, 
$D(r_2)=(0,30)$, $P(r_2)=\text{\normalfont uniform}(0,30)$, 
$D(r_3)=(0,15)$, $P(r_3)=\text{\normalfont uniform}(0,15)$, 
$D(r_4)=(0,20)$, $P(r_2)=\text{\normalfont uniform}(0,20)$, 
$C=\{c_1:C_1\geq X_1 r_1 + X_2 r_2, c_2:C_2\geq X_1 r_3 + X_2 r_4\}$, $\beta_{c_1}=\beta_{c_2}=0.7$, and $L=[\langle V,S\rangle]$.
$C_1=245$ and $C_2=215$ are constants.
We set $\alpha=0.9$ and $\vartheta=0.05$. We computed analytically the true boundaries of $c_1$ and $c_2$ (see \cite{RePEc:spr:stpapr:v:50:y:2009:i:1:p:171-175,citeulike:11042202}), each of which is denoted by a dashed line in Fig. \ref{fig:graph_region_2848} and \ref{fig:graph_region_348}. We also computed confidence bands around these two dashed lines. The upper confidence band is the set of solutions that provide a satisfaction probability of exactly $\beta_i-\vartheta$; the lower confidence band is the set of solutions that provide a satisfaction probability of exactly $\beta_i+\vartheta$. 
We apply Definition \ref{sample_size_alpha_theta_solution_set} to compute the number of samples $N=2848$ required to obtain an $(\alpha,\vartheta)$-solution set to $\mathcal{P}$, which is shown in Fig. \ref{fig:graph_region_2848}. 
\end{example}

\subsection{Approximating ($\alpha$,$\vartheta$)-solution sets}
\label{sec:alpha_theta_solution_set_approx}

Bonferroni's correction is known to be conservative. In particular, as we have seen, this correction assumes events $\tau_i=1$ and $\tau_j=1$ are mutually exclusive for all $i,j=1,\ldots,T$, $i\neq j$. In other words, we are assuming that assignment misclassifications are mutually exclusive. In practice, in an SCSP sets of assignments are often misclassified together depending on random variables realisations. For this reason a correction such as the one introduced in Definition \ref{sample_size_alpha_theta_solution_set} will generally be too conservative and will lead to a sample size much larger than the one strictly needed to obtain an $(\alpha,\vartheta)$-solution set. This fact is well known in statistics and a number of adjusted corrections have been proposed to account for correlated errors \citep[see e.g.][]{lehmann2005,citeulike:13411581}. 

In what follows, we shall therefore adopt a less conservative approximate correction strategy. To the best of our knowledge no similar correction has been discussed in the literature. In our computational study (Section \ref{sec:comp}) we will demonstrate the effectiveness of this technique. Of course, the investigation of other less conservative and possibly exact correction strategies, ideally borrowed from established results in statistics, is an interesting direction for future research.

Consider the general case in which constraint $h$ constrains all $m$ random variables in $S$.
\begin{lem}\label{lem:det_test}
Given realisations $\{\bar{s}^1_1,\ldots,\bar{s}^1_{m}\}$, $\{
\bar{s}^2_1,\ldots,\bar{s}^2_{m}\}$, $\ldots$,
 $\{\bar{s}^N_1,\ldots,\bar{s}^N_{m}\}$,
where $\bar{s}^k_j$ is the realised value for random variable $j$ observed in the $k$-th 
set of realisations, $\tau_i$ is a deterministic test.
\end{lem}

\begin{prop}\label{prop:random_function}
$\tau_i$ is a function of random variables $s_1,\ldots,s_{m}$ and of $N$. 
\end{prop}
\begin{proof}
Immediately follows from Lemma \ref{lem:det_test} and from the fact that $s^1_j,\ldots,s^N_j$ are $N$ i.i.d. random variables.
\end{proof}

\begin{prop}
The probability that at least one of $\tau_1,\ldots,\tau_T$ takes value 1 is uniquely determined by the probability distributions of $s_1,\ldots,s_{m}$ and the number of samples $N$.
\end{prop}
\begin{proof}
Follows from the definition of $\tau$, Lemma \ref{lem:det_test} and Proposition \ref{prop:random_function}.
\end{proof}

In practice, this means that assignment misclassifications in a sampled SCSP, e.g. events $\tau_i=1$ and $\tau_j=1$, depend on realisations of one or more random variables $s_1,\ldots,s_{m}$; note that $m$ is generally much smaller than $T$. Since the multivariate random variable $\{\tau_1,\ldots,\tau_T\}$ is a (deterministic) function of the multivariate random variable $\{s_1,\ldots,s_{m}\}$ and of $N$ (Proposition \ref{prop:random_function}), and since in the previous section we have fully characterised the marginal probability distribution of a test $\tau_i$, the probability that $\tau_1,\ldots,\tau_T$ are all equal to 0 is approximately bounded from below by $1 - m(1-\alpha)$; i.e. we correct for at most $m$ mutually exclusive misclassifications induced by random variables $s_1,\ldots,s_{m}$ and we assume that all remaining misclassifications depend on one or more of these. Once more we introduce a correction for each chance constraint $h\in C$. Let $m_h$ be the number of random variables constrained by $h$, to compute an approximate $(\alpha,\vartheta)$-solution set we introduce the following correction while computing $N$.

\begin{definition}\label{sample_size_alpha_theta_solution_set_approx}
$N$ is computed as the minimum value for which \[\max(p^{\beta}_{\text{ub}}-\beta,\beta-p^{\beta}_{\text{lb}})\leq \vartheta,\]
where $p^{\beta}_{\text{lb}}$ and $p^{\beta}_{\text{ub}}$ are the single-sided Clopper-Pearson confidence interval bounds for a confidence probability $\widehat{\alpha}$, where  \[\widehat{\alpha}= 1-\frac{1-\alpha}{\sum_{h\in C}m_h},\] and
$\mbox{\it round}(\beta N)$ ``successes'' in $N$ trials.
\end{definition}

A set of policy trees that are solutions to $\widehat{\mathcal{P}}$ for a sample size $N$ computed as discussed in Definition \ref{sample_size_alpha_theta_solution_set_approx} is an approximate $(\alpha,\vartheta)$-solution set for $\mathcal{P}$.

\begin{example}\label{ex:example_6}
Consider once more the SCSP in Example \ref{ex:example_5}. We apply Definition \ref{sample_size_alpha_theta_solution_set_approx} to compute the number of samples $N=348$ required to obtain an approximate $(\alpha,\vartheta)$-solution set to $\mathcal{P}$, which is shown in Fig. \ref{fig:graph_region_348}; note that there are two constraints each of which constrains two random variables. 
To assess the quality of this approximation, we generated 1000 different instances and analytically inspected, for each of them, if the $(\alpha,\vartheta)$-solution set generated was fully contained within the upper confidence band in Fig. \ref{fig:graph_region_348}; the result of this simulation study revealed that the $(\alpha,\vartheta)$-solution set was not fully contained within the upper confidence band with probability 0.894, 0.95 confidence interval $(0.873,0.912)$; this misclassification rate is in line with the prescribed $\alpha$. Finally, it is worth noting that the random boundary of an $(\alpha,\vartheta)$-solution set remains within the channel identified by the two solid confidence bands with probability at least $1-2(1-\alpha)$. 
\end{example}

\section{Stochastic constraint optimisation problems}\label{sec:scop}

The concepts introduced in sections \ref{sec:sampled_scsps_solutions} and \ref{sec:alpha_theta_solution_set}  can be employed to approximate optimal solutions to sampled SCOPs.
In this setting, we must distinguish two possible cases: the case in which the objective function is deterministic and that in which the objective function is stochastic.
If the objective function is deterministic, it is possible to exploit the results in section \ref{sec:alpha_theta_solution_set} to obtain a confidence interval for the cost/profit of an optimal plan. 
Without loss of generality, we discuss the case in which our aim is to maximise a deterministic objective function $f$ of the decision variables in $V$. Consider an SCOP $\mathcal{P}=\langle V,S,D,P,C,\beta_c, L, f \rangle$. Choose $\alpha$ and $\vartheta$ and construct two new SCOPs: $\mathcal{P}_{\text{lb}}=\langle V,S,D,P,C,\beta^1_c, L, f \rangle$, where 
for all $c\in C$, $\beta^1_c=\beta_c+\vartheta$; and $\mathcal{P}_{\text{ub}}=\langle V,S,D,P,C,\beta^2_c, L, f \rangle$, where for all $c\in C$, $\beta^2_c=\beta_c-\vartheta$.
\begin{prop}\label{prop:bounds}
An $(\alpha,\vartheta)$-solution set to $\mathcal{P}_{\text{lb}}$ underestimates the true optimal profit with probability greater or equal to $\alpha$; an $(\alpha,\vartheta)$-solution set to $\mathcal{P}_{\text{ub}}$ overestimates the true optimal profit with probability greater or equal to $\alpha$.
\end{prop}
\begin{proof}
The proof follows from Definition \ref{alpha_theta_solution_set}.
\end{proof}
Proposition \ref{prop:bounds} can be exploited to generate a confidence interval for the true optimal profit via a binomial reasoning. We solve $M$ independently generated instances of $\mathcal{P}_{\text{lb}}$ and store the optimal profit obtained for each of these instances into an array $K_{\text{lb}}$ sorted in ascending order; we solve $M$ independently generated instances of $\mathcal{P}_{\text{ub}}$ and store the optimal profit obtained for each of these instances into an array $K_{\text{ub}}$ sorted in ascending order. Let $\text{bin}^{-1}(M,\alpha)$ be the inverse cumulative distribution of a binomial distribution with $M$ trials and a success probability $\alpha$; let $k_{\text{lb}}$ be the $(1-\alpha)/2$-quantile of this distribution; finally, let $k_{\text{ub}}$ be the $1-(1-\alpha)/2$-quantile of $\text{bin}^{-1}(M,1-\alpha)$. With confidence $\alpha$ element at position $k_{\text{lb}}$ of $K_{\text{lb}}$ is a lower bound and element at position $k_{\text{ub}}+1$ of $K_{\text{ub}}$  is an upper bound to the true optimal cost.\footnote{
Elements of $K_i$ are indexed as follows: $1,\ldots,|K_i|$. Note that in statistics the $k^{\text{th}}$-smallest value of a statistical sample is known as $k^{\text{th}}$ order statistic \cite{citeulike:2594180}.} 

\begin{example}\label{ex:example_4a}
We transform the SCSP in Example \ref{ex:example_5} into two SCOPs $\mathcal{P}_{\text{lb}}$ and $\mathcal{P}_{\text{ub}}$ that maximise the objective function $f(X_1,X_2)=X_1+2X_2$. In other words, we assume the profit per unit of $X_1$ is $1$ and the profit per unit of $X_2$ is $2$. By choosing $M=20$ we obtain the $\alpha$ confidence interval $(282,304)$ for the true optimal profit $293$; if we reduce $\vartheta$ to 0.01 the interval shrinks considerably to $(290,295)$.
\end{example}

If the objective function is stochastic there is no unique way to proceed. For instance, based on the available samples one may derive standard confidence intervals for the expected value of a stochastic expression based on the Student's $t$ distribution and then compare solutions or partial assignments by comparing upper or lower limits of these intervals. The decision maker must of course choose a suitable confidence level $\alpha$ associated with this estimation. An example of a filtering algorithm that may be employed in such context is discussed in \ref{sec:appendix_I}. This algorithm is designed to handle the situation in which the objective is to minimise/maximise the expected value of some expression involving decision and random variables. Different algorithms must be designed if the objective involves a different operator, e.g. variance. Our algorithm distinguishes the case in which we are trying to determine an upper or a lower bound for the expected cost of an optimal solution. It then exploits the sampling distribution (i.e. Student's $t$ distribution) of the expected total profit/cost and filters values based on upper/lower confidence limits obtained via this distribution. For instance, if our aim is to determine an upper bound for the optimal profit (problem type $\mathcal{P}_{\text{ub}}$), our algorithm will simply compare the upper confidence limits of the expected profit of two assignments and retain the assignment with the highest upper confidence limit. We will make use of this propagator to solve the models discussed in section \ref{sec:comp}. 

Finally, one should note that an alternative strategy may instead compare not only the upper confidence limits, but the whole intervals. An assignment would then provide a lower/higher profit than another if and only if their profit confidence intervals do not overlap. However, due to the complexity of the filtering logic that would be required in this case, we prefer to leave this discussion as future work. 

\section{Connections with statistics}\label{sec:connections}

To better understand the concepts just introduced, it is worth discussing the connection between
the approach introduced and hypothesis testing in statistical analysis. Let us assume that 
our {\em null hypothesis} ($H_0$), in statistical sense, is that an assignment is feasible. 
According to classical hypothesis testing we may have four cases, as illustrated in Table \ref{tab:errors}.
We may have a feasible assignment at hand ($H_0$ true) and we may incorrectly filter it (Type I error);
or we may be operating on an infeasible assignment ($H_0$ false) and we may fail to reject it (Type II error).

\begin{table}[htdp]
\begin{center}
\begin{tabular}{c|cc}
&$H_0$ is true&$H_0$ is false\\
\hline
Reject $H_0$&Type I error&Correct outcome\\
			&(false positive)& (true positive)\\
Fail to reject $H_0$&Correct outcome &Type II error \\
			&(true negative)&(false negative)
\end{tabular}
\end{center}
\caption{Type I and Type II errors in statistics}
\label{tab:errors}
\end{table}%

In clinical trials or quality control, it is key to control the rate of Type I errors. It is undesirable to put under treatment a healthy a patient or to discard a functioning expensive machine. However, there are cases in which controlling Type II errors is essential.  For example, aerospace engineers would prefer to scrap a functioning electronic circuit  
than to use one that is actually broken on a spacecraft;
in such a situation a Type I error raises the budget, but a Type II error would put at risk the entire mission. 
In general, minimising Type I and Type II errors is not a simple matter. If one tries to reduce the rate of occurrence for Type I errors, the direct consequence is typically an increase in the observed rate for Type II errors and vice-versa. So in practice, one tries to control either Type I or Type II errors and, if the rate of the type that is not controlled is too high, then one increases the sample size. 

In our specific case it is clearly essential to control the rate of Type II errors, which are more delicate 
than Type I errors. Making a Type II error means retaining an infeasible assignment, which is what we 
want to avoid as much as possible. Making a Type I error means discarding a feasible solution, which may impact optimality for an optimisation problem, or may lead to an empty solution space. Since our approach is essentially a heuristic, it is clear that both these issues --- a poor solution quality or an empty solution space --- are acceptable and should be dealt with by increasing the number of samples. 

\section{Computational experience}\label{sec:comp}

The aim of this section is to provide numerical insights on the theoretical framework introduced and particularly on the concept of ($\alpha$,$\vartheta$)-solution set and on its applications to find approximate solution to SCSPs and SCOPs.
In our numerical study we will consider three well-known problems: the static stochastic knapsack (Section \ref{sec:ssmkp}), the stochastic multiprocessor scheduling problem with release time and deadlines (Section \ref{sec:smps}), and the static stochastic lot-sizing problem (Section \ref{sec:sslsp}). The first and the third problems are single-stage, while the second is two-stage.
In Section \ref{sec:feasibility} we will generate approximate ($\alpha$,$\vartheta$)-solution sets using Definition \ref{sample_size_alpha_theta_solution_set_approx} for the first two problems and show numerically that, with probability greater than or equal to $\alpha$, the approach we discussed generates solution sets that satisfy chance constraints in the model with a margin of error $\vartheta$. 
In Section \ref{sec:optimality} we will numerically illustrate that the upper and lower profit/cost bounds obtained with the approach outlined in Section \ref{sec:scop} comply with the prescribed confidence level $\alpha$. We will also show the behaviour of the optimality gap as a function of the chosen error threshold $\vartheta$ and number $M$ of independently generated instances of $\mathcal{P}_{\text{lb}}$ and $\mathcal{P}_{\text{ub}}$. Finally, in Section \ref{sec:comp_efficiency} we will investigate computational efficiency and scalability. All our experiments were performed by using Choco \cite{ocre} on an Intel Xeon(r) CPU @ 3.50 Ghz with 16GB of RAM. 

\subsection{Static stochastic knapsack}\label{sec:ssmkp}

The knapsack problem \cite{citeulike:13259752} is a well-known combinatorial optimisation problem. The decision maker is given a set of objects each of which is associated with a weight and a profit. The aim is then to select a subset of these objects that fit into a given capacity and bring the maximum profit. There are several possible stochastic variants of the knapsack problem. Stochastic versions of the knapsack problem can be classified into static or dynamic. In the static stochastic knapsack problem, see e.g. \cite{ksh01}, object weights and/or profits are random and the decision maker must choose, before observing any of their weights/profits, a subset of these objects that maximises a given objective, e.g. the expected profit, while meeting a restriction, e.g. a chance constraint, on the given capacity. Conversely, in the dynamic stochastic knapsack, see e.g. \cite{Kleywegt1998}, the decision maker selects an object and immediately observes its weight and/or profit; based on this information she can then decide whether to select or not other objects.

In our computational study we will consider the SCSP presented in Fig. \ref{model:ssmkp}, i.e. a static stochastic multiple knapsack (SSMKP). In this problem we have a set of $N$ types of objects; there are $D$ objects of type $i$ available. Each object of type $i$ is associated random ``coefficients'' $s^k_i$ that appear in the context of $G$ chance constraints --- this set of coefficients is generally denoted as {\em stochastic technology matrix} \cite{citeulike:8219880}; without loss of generality, these coefficients follow a Poisson distribution with mean $\lambda^k_i$.\footnote{For a discussion on statistic stochastic knapsack problems with Poisson resource requirements see e.g. \cite{citeulike:5128004}.} The first $L$ of these chance constraints are of type (1), i.e. they can be seen as ``capacity restrictions'' with respect to a target capacity $C_k$, and they should be satisfied with probability $\beta$. In the context of the first $L$ chance constraints $s^k_i$ represents the ``weight'' of item $i$ in chance constraint $k$. The remaining $G-L$ chance constraints are of type (2), i.e. they can be seen as ``minimum production requirements'' with respect to a target level $C_k$, and again they should be satisfied with probability $\beta$. In the context of the remaining $G-L$ chance constraints $s^k_i$ represents the ``production contribution'' of item $i$ in chance constraint $k$.

Our aim is to determine the feasible region of the problem, i.e. the set of assignments that satisfy constraints (1) and (2).

\begin{figure}[t]
\begin{center}
\framebox{
        \begin{tabular}{ll}
        \mbox{{\bf Constraints:}} \\
        (1) $~\Pr\{s^k_1x_1+\ldots+s^k_Nx_N\leq C^k\} \geq \beta$&$k=1,\ldots,L$\\
        (2) $~\Pr\{s^k_1x_1+\ldots+s^k_Nx_N\geq C^k\} \geq \beta$&$k=L+1,\ldots,G$\\\\
        \mbox{{\bf Decision variables:}} \\
        $~~~~x_i \in \{0,\ldots,D\}$&$i=1,\ldots,N$\\\\
        \mbox{{\bf Random variables:}} \\
        $~~~~s^k_i \leftarrow \text{Poisson}(\lambda^k_i)$&$i=1,\ldots,N; k=1,\ldots,G$ \\\\
        \mbox{{\bf Stage structure:}} \\
        $~~~~V_1=\{x_1,\ldots,x_N\}$\\
        $~~~~S_1=\{s^1_1,\ldots,s^k_n,\ldots,s^G_N\}$\\
        $~~~~L=[\langle V_1,S_1 \rangle]$
        \end{tabular}
        }
        \end{center}
       \caption{The static stochastic multiple knapsack as an SCSP} 
       \label{model:ssmkp}  
\end{figure}

\begin{figure}[t]
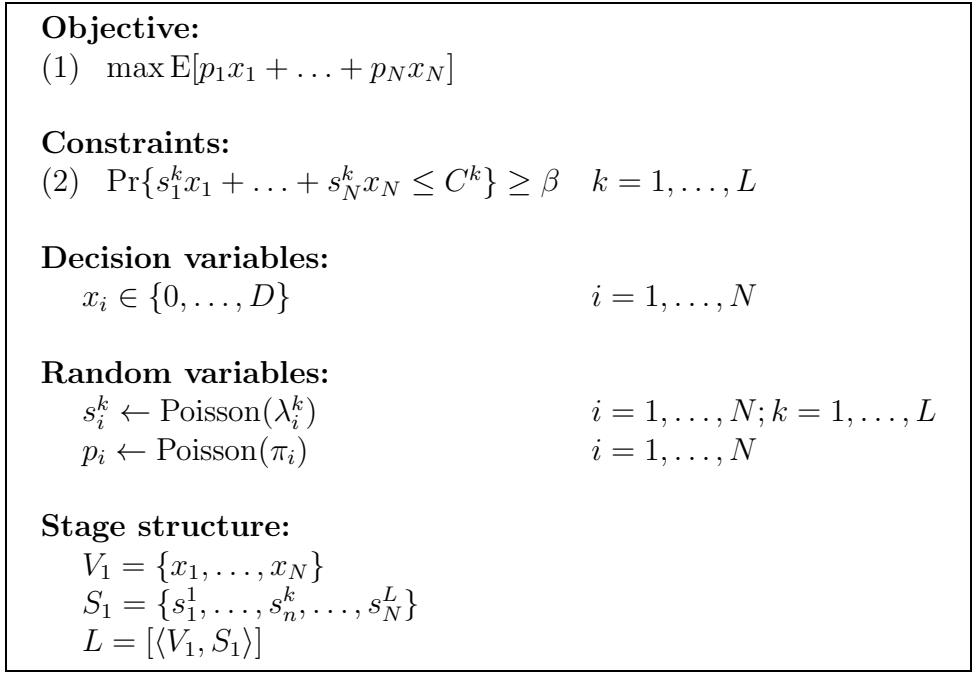

\begin{center}
\framebox{
        \begin{tabular}{ll}
        \mbox{{\bf Objective:}} \\
        (1) $~\max \text{E}[p_1x_1+\ldots+p_Nx_N]$\\\\
        \mbox{{\bf Constraints:}} \\
        (2) $~\Pr\{s^k_1x_1+\ldots+s^k_Nx_N\leq C^k\} \geq \beta$&$k=1,\ldots,L$\\\\
        \mbox{{\bf Decision variables:}} \\
        $~~~~x_i \in \{0,\ldots,D\}$&$i=1,\ldots,N$\\\\
        \mbox{{\bf Random variables:}} \\
        $~~~~s^k_i \leftarrow \text{Poisson}(\lambda^k_i)$&$i=1,\ldots,N; k=1,\ldots,L$ \\
        $~~~~p_i \leftarrow \text{Poisson}(\pi_i)$&$i=1,\ldots,N$ \\\\        
        \mbox{{\bf Stage structure:}} \\
        $~~~~V_1=\{x_1,\ldots,x_N\}$\\
        $~~~~S_1=\{s^1_1,\ldots,s^k_n,\ldots,s^L_N\}$\\
        $~~~~L=[\langle V_1,S_1 \rangle]$
        \end{tabular}
        }
        \end{center}
       \caption{The static stochastic multiple knapsack as an SCOP} 
       \label{model:ssmkp_opt}
\end{figure}

We will also consider an optimisation version of the problem (Fig. \ref{model:ssmkp_opt}) in which our aim is to determine what subset of the $N$ objects in the problem maximises the expected total profit while satisfying all chance constraints. For each object $i$, we therefore introduce a random ``profit'' $p_i$, which follows a Poisson distribution with mean $\pi_i$; once more the choice of the distribution is made without loss of generality.

\subsection{Stochastic multiprocessor scheduling problem with release time and deadlines}\label{sec:smps}

We consider a multiprocessor scheduling problem (MPSP, see \cite{citeulike:574128}, p. 238). The problem consists in finding a feasible schedule to process a set of $K$ orders (or jobs) using $m$ processors, where $m\leq P$. Processing an order $k$ can only begin after the release date $r_k$ and must be completed at the latest by the due date $d_k$. Order $k$ requires a certain capacity $c_k$ --- expressed in terms of the number of processors --- to be processed. The processing time of order $k$ is $t_k$. The problem just described is well known in scheduling and it is fully deterministic and can easily and compactly be modelled using the \texttt{cumulative} constraint \cite{citeulike:8385205}. Let the height of a task $k$ be $c_k$. This constraint considers a set of tasks and enforces that at each point in time the total height of the set of tasks that overlap that point does not exceed a given limit $m$. A task $k$ overlaps a point $i$ if and only if its origin $s_k$ is less than or equal to $i$, and its end $e_k$ is strictly greater than $i$. This constraint also imposes, for each task $k$, the constraint $s_k$+$t_k$=$e_k$.

However, in reality, some parameters of this problem are uncertain in nature.  Jobs may take longer than expected, some processors may break down and become unavailable, the release and due dates may be delayed, etc. 
In order to better model this problem a number of stochastic generalizations may be considered such as uncertain release date $r_k$; uncertain due date $d_k$; uncertain processing capacity $c_k$; uncertain processing time $t_k$; and uncertain number $m$ of available processors; and every possible combination stemming from these cases. 

We will consider the following stochastic constraint programming
formulation of the stochastic multiprocessor scheduling problem (SMPSP), in which only processing time
$t_k$ for order $k$ is uncertain; this is shown in Fig. \ref{model:ss_scp_multiprocessor_scheduling}.
\begin{figure}[t]
\begin{center}
\framebox{
        \begin{tabular}{ll}
        \mbox{{\bf Constraints:}}\\
        $~~~~(1)~~\Pr\left\{\mathrm{cumulative}(s,e,t,c,m) \right\} \geq \beta$ \\\\
        \mbox{{\bf Decision variables:}} \\
        ~~~~$s_k \in \{r_k,\ldots,d_k\}$,&$\forall k \in 1,\ldots,K$\\
        ~~~~$e_k \in \{r_k,\ldots,d_k\}$,&$\forall k \in 1,\ldots,K$\\\\
        \mbox{{\bf Stochastic variables:}} \\
        ~~~~$t_k \rightarrow$ \text{Poisson}$(\lambda_k)$&$\forall k \in 1,\ldots,K$\\\\
        \mbox{{\bf Stage structure:}} \\
        ~~~~$V_1=\{s_1,s_2,\ldots,s_K\}$~~~$S_1=\{t_1,t_2\ldots,t_K\}$\\
        ~~~~$V_2=\{e_1,e_2,\ldots,e_K\}$~~~$S_2=\{\}$\\
        ~~~~$L=[\langle V_1,S_1\rangle,\langle V_2,S_2\rangle]$
        \end{tabular}
        }
\end{center}
       \caption{An SCSP for the stochastic multiprocessor scheduling problem with release time and deadlines}
       \label{model:ss_scp_multiprocessor_scheduling}
\end{figure}
In this model, decision variables  $s_k$ and $e_k$ denote the start
time and the completion time of each job $k$, respectively. The processing time $t_k$ of each job $k$ is modeled as a Poisson distributed random variable with mean $\lambda_k$. In contrast to the problem presented in Section \ref{sec:ssmkp}, this model is a two-stage SCSP. In the first stage, we decide on the start time of each job then we observe the realisation of the processing time. In the second stage the completion times are decided.  Under this stage structure, constraint (1) enforces that the \emph{probability} of not exceeding the given deadline for each job and the number of available processors $m$ stays above the specified threshold $\beta$. More specifically, this constraint is a global chance constraint embedding a well-known global constraint: the \texttt{cumulative} constraint
\cite{citeulike:8385205}. This constraint can be filtered using the general purpose method discussed in \cite{citeulike:10689458} .

In our computational study we will also consider an optimisation version of the above problem in which we aim to minimise the latest start time.

\subsection{Static stochastic lot-sizing}\label{sec:sslsp}

The last problem we will consider in our computational study is the single-item stochastic lot-sizing problem introduced in \cite{bt88}. A SCOP for this problem is shown in in Fig. \ref{model:lotsizing}. The decision maker faces a finite horizon of $T$ periods and a random demand $d_t$ in each period; which, without loss of generality, we will consider Poisson distributed with mean $\lambda_t$. 
There is a fixed cost $a$ for placing an order of size $0<Q_t\leq C$ in period $t$. An order placed in period $t$ is delivered immediately at the beginning of the period, before demand occurs. Binary decision variable $\delta_t$ is set to zero if no order is placed (3). There is a holding cost $h$ charged on items that are carried over from one period to the next. Finally, the decision maker must comply with a service level restriction (2) stating that the net inventory at the end of each period should be nonnegative with probability at least $\beta$. The aim is to meet these service level restrictions while minimising the expected total cost (1). 

The authors in \cite{bt88} describe a range of control policies that can be used to control such a system. In our study, we will adopt the static uncertainty policy, which fixes all $Q_t$ and $\delta_t$ at the beginning of the planning horizon, before demand is observed. Note that other strategies discussed in \cite{bt88}, i.e. dynamic uncertainty and static-dynamic uncertainty, can be easily captured by modifying the stage structure of the SCOP. In what follows, we shall refer to this problem as the static stochastic lot-sizing problem (SSLSP).
\begin{figure}[t]
\begin{center}
\framebox{
        \begin{tabular}{ll}
        \mbox{{\bf Objective:}} \\
        (1) $~\min \text{E}[\sum_{t=1}^N(a\delta_t+h\sum_{j=1}^t(Q_t-d_t))]$\\\\
        \mbox{{\bf Constraints:}} \\
        (2) $~\Pr\{\sum_{j=1}^t(Q_t-d_t)\geq0\} \geq \beta$&$t=1,\ldots,T$\\
        (3) $~\delta_t = 0 \implies Q_t = 0$&$t=1,\ldots,T$\\\\
        \mbox{{\bf Decision variables:}} \\
        $~~~~\delta_t \in \{0,1\}$&$t=1,\ldots,T$\\
        $~~~~Q_t \in \{0,\ldots,C\}$&$t=1,\ldots,T$\\\\
        \mbox{{\bf Random variables:}} \\
        $~~~~d_t \leftarrow \text{Poisson}(\lambda_t)$&$t=1,\ldots,T$ \\\\
        \mbox{{\bf Stage structure:}} \\
        $~~~~V_1=\{Q_1,\ldots,Q_T,\delta_1,\ldots,\delta_T\}$\\
        $~~~~S_1=\{d_1,\ldots,d_T\}$\\
        $~~~~L=[\langle V_1,S_1 \rangle]$
        \end{tabular}
        }
\end{center}
      \caption{The stochastic lot-sizing problem in \cite{bt88} as an SCOP (static uncertainty strategy)} 	     
       \label{model:lotsizing}
\end{figure}

\subsection{Feasibility}\label{sec:feasibility}

In Section \ref{sec:alpha_theta_solution_set} we introduced the notion of approximate ($\alpha$,$\vartheta$)-solution set. We will now present a computational analysis for the SCSPs presented in Sections \ref{sec:ssmkp} and \ref{sec:smps} demonstrating that, with probability $\alpha$, our approach generates solution sets that satisfy chance constraints in the model with a margin of error $\vartheta$. 

We considered thirty randomly generated small instances of the single stage problem in Fig. \ref{model:ssmkp} in which $N=2$, $L=2$, $G=3$, $D=250$ and $\beta=0.7$. Means $\lambda_i^k$ of random variables in the model were integer numbers uniformly distributed between 10 and 20 for constraints (1) and between 20 and 30 for constraints (2). Right hand side constants $C^k$ were integer numbers uniformly distributed between 1500 and 2000 for constraints (1) and between 2500 and 3000 for constraints (2). We fixed $\alpha=0.9$ and $\vartheta=0.2$; according to Definition \ref{sample_size_alpha_theta_solution_set_approx} this choice led to a sample size of 31.

We also considered thirty randomly generated small instances of the two stage problem in Fig. \ref{model:ss_scp_multiprocessor_scheduling} in which $K=2$ and $\beta=0.6$; $r_k$ and $d_k$, which represent job $k$ release time and deadline, were set to $0$ and $4$, respectively. Capacity requirements $c_k$ were generated as integer numbers uniformly distributed between 1 and 2. Finally, expected task durations $\lambda_k$ were generated as uniformly distributed numbers between 1 and 3; the maximum number of processors $P$ was set to 3. We fixed $\alpha=0.9$ and $\vartheta=0.35$, this choice led to a sample size of 6.

Instances were small since in our analysis we generated the complete set of feasible assignments of the respective sampled SCSP, i.e. an $(\alpha,\vartheta)$-solution set, which for the two-stage problem in Fig. \ref{model:ss_scp_multiprocessor_scheduling} was generally extremely large, in the order of tens of thousands of solutions. Feasibility of each of these assignment with respect to the original SCSP was then assessed via Monte Carlo simulation; the number of simulation runs was set in such a way as to guarantee a margin of error of $\vartheta/10$ with a confidence level of $0.9$ --- so that the Monte Carlo simulation error is an order of magnitude smaller than the approximation error associated with the ($\alpha$,$\vartheta$)-solution set obtained. 

\begin{figure}[ht]
\begin{center}
\[\begin{array}{cc}
\resizebox{0.48\columnwidth}{!}{
\begin{tikzpicture}
\begin{axis}[
title={SSMKP},
xtick = {1,10,20,30},
xlabel=Instance,
ylabel=Frequency]

\addplot[smooth,densely dotted,mark=none,
domain=-1:31,samples=40]
{0.9};

\addplot+[black,
mark=-,
only marks,
error bars/.cd,
y dir=both,y explicit]
coordinates {
(	1	,	1	)	+-	(0,	0.003682084	)
(	2	,	0.943	)	+-	(0,	0.016222729	)
(	3	,	1	)	+-	(0,	0.003682084	)
(	4	,	1	)	+-	(0,	0.003682084	)
(	5	,	1	)	+-	(0,	0.003682084	)
(	6	,	1	)	+-	(0,	0.003682084	)
(	7	,	0.954	)	+-	(0,	0.014883553	)
(	8	,	1	)	+-	(0,	0.003682084	)
(	9	,	1	)	+-	(0,	0.003682084	)
(	10	,	0.914	)	+-	(0,	0.019120268	)
(	11	,	0.966	)	+-	(0,	0.013189439	)
(	12	,	0.878	)	+-	(0,	0.021904219	)
(	13	,	0.982	)	+-	(0,	0.010299041	)
(	14	,	0.94	)	+-	(0,	0.01656049	)
(	15	,	1	)	+-	(0,	0.003682084	)
(	16	,	0.903	)	+-	(0,	0.020047844	)
(	17	,	0.988	)	+-	(0,	0.00886768	)
(	18	,	0.891	)	+-	(0,	0.020978365	)
(	19	,	0.989	)	+-	(0,	0.008596628	)
(	20	,	1	)	+-	(0,	0.003682084	)
(	21	,	0.964	)	+-	(0,	0.013492912	)
(	22	,	0.991	)	+-	(0,	0.008015783	)
(	23	,	1	)	+-	(0,	0.003682084	)
(	24	,	0.915	)	+-	(0,	0.019031986	)
(	25	,	1	)	+-	(0,	0.003682084	)
(	26	,	0.921	)	+-	(0,	0.018487095	)
(	27	,	0.884	)	+-	(0,	0.021486751	)
(	28	,	0.995	)	+-	(0,	0.006629471	)
(	29	,	0.972	)	+-	(0,	0.012214887	)
(	30	,	0.9	)	+-	(0,	0.020287937	)
};
\end{axis}
\end{tikzpicture}
}
&
\resizebox{0.48\columnwidth}{!}{
\begin{tikzpicture}
\begin{axis}[
title={SMPSP},
xtick = {1,10,20,30},
xlabel=Instance,
ylabel=Frequency]

\addplot[smooth,densely dotted,mark=none,
domain=-1:31,samples=40]
{0.9};

\addplot+[black,
mark=-,
only marks,
error bars/.cd,
y dir=both,y explicit]
coordinates {
(	1	,	0.953	)	+-	(0,	0.015012631	)
(	2	,	0.954	)	+-	(0,	0.014883553	)
(	3	,	0.972	)	+-	(0,	0.012214887	)
(	4	,	0.996	)	+-	(0,	0.006209665	)
(	5	,	0.982	)	+-	(0,	0.010299041	)
(	6	,	0.956	)	+-	(0,	0.014620419	)
(	7	,	0.969	)	+-	(0,	0.012715085	)
(	8	,	0.947	)	+-	(0,	0.015755183	)
(	9	,	0.997	)	+-	(0,	0.005742023	)
(	10	,	0.978	)	+-	(0,	0.011119665	)
(	11	,	0.991	)	+-	(0,	0.008015783	)
(	12	,	0.944	)	+-	(0,	0.016107746	)
(	13	,	0.983	)	+-	(0,	0.010079497	)
(	14	,	0.951	)	+-	(0,	0.015266053	)
(	15	,	0.979	)	+-	(0,	0.010922335	)
(	16	,	0.999	)	+-	(0,	0.004558924	)
(	17	,	0.973	)	+-	(0,	0.012041747	)
(	18	,	0.979	)	+-	(0,	0.010922335	)
(	19	,	0.997	)	+-	(0,	0.005742023	)
(	20	,	0.96	)	+-	(0,	0.014072697	)
(	21	,	0.971	)	+-	(0,	0.012384719	)
(	22	,	0.976	)	+-	(0,	0.011500639	)
(	23	,	0.97	)	+-	(0,	0.012551402	)
(	24	,	0.979	)	+-	(0,	0.010922335	)
(	25	,	0.989	)	+-	(0,	0.008596628	)
(	26	,	0.971	)	+-	(0,	0.012384719	)
(	27	,	0.964	)	+-	(0,	0.013492912	)
(	28	,	0.955	)	+-	(0,	0.014752833	)
(	29	,	0.984	)	+-	(0,	0.009853249	)
(	30	,	0.977	)	+-	(0,	0.011312338	)
};
\end{axis}
\end{tikzpicture}
}
\end{array}\]
\end{center}
\caption{Frequency of event ``all feasible assignments of the sampled SCSP are feasible with respect to the original SCSP within the given tolerance threshold $\vartheta$'' over 1000 sampled SCSPs; together with the frequency, we report the associated confidence interval (confidence level of $0.95$)}
\label{fig:alpha_theta_solutions}
\end{figure}
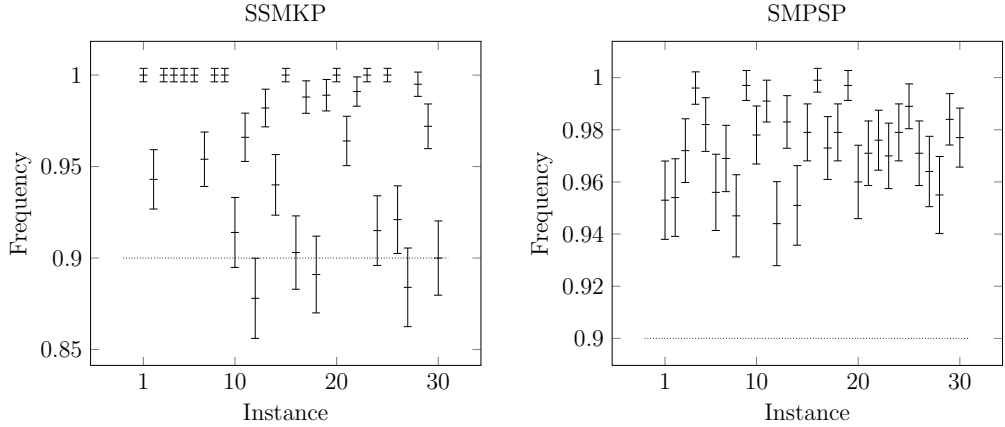

To numerically investigate if those computed are effectively $(\alpha,\vartheta)$-solution sets, for each of the above sixty instances, we repeatedly solved 1000 sampled SCSPs and computed the frequency of event $e$: ``{\bf all feasible assignments} of the sampled SCSP are feasible with respect to the original SCSP within the given tolerance threshold $\vartheta$.'' In Fig. \ref{fig:alpha_theta_solutions}, for both problems and for each instance, we report the frequency of event $e$ and the associated confidence intervals (confidence level of $0.95$). These frequencies, are in line with the claim that those computed are $(\alpha,\vartheta)$-solution sets, for the given $\alpha=0.9$. Note that our aim is to control Type-II errors (an infeasible assignment regarded as feasible), and not Type-I errors (a discarded and yet feasible assignment); for this reason if the sampled SCSP admitted no solution, this was regarded as a degenerate case in which all feasible assignments (i.e. none) of the sampled SCSP were feasible with respect to the original SCSP within the given tolerance threshold $\vartheta$. Finally, it is worth observing that some of the frequencies observed in Fig. \ref{fig:alpha_theta_solutions} are strictly greater than the prescribed value $\alpha$. This is due to the fact that only assignments providing a satisfaction probability of exactly $\beta-\vartheta$ are correctly classified as infeasible with probability $\alpha$. However, given the discrete nature of the assignment space, it is likely that instances may not feature any such assignment. Assignments providing a satisfaction probability strictly less than $\beta-\vartheta$ are correctly classified as infeasible with probability strictly greater than $\alpha$. In addition to this, when a model features multiple chance constraints, 
Bonferroni's correction, which is free of correlation and distribution assumptions, might generate a conservative --- i.e. strictly larger than needed --- sample size.

\subsection{Optimality}\label{sec:optimality}

We considered fifty randomly generated small instances of the problem in Fig. \ref{model:ssmkp_opt} (SSMKP) in which $N=10$, $L=2$, $D=1$ and $\beta=0.9$. Means $\lambda_i^k$ of random variables in the model were integer numbers uniformly distributed between 10 and 20 for constraints (1). Right hand side constants $C^k$ were integer numbers uniformly distributed between 100 and 200 for constraints (1). Means $\pi_i$ were all set to 10. 

We also considered fifty randomly generated small instances of the problem in Fig. \ref{model:lotsizing} (SSLSP) in which $T=5$, $h=1$, $a=10$, $C=100$, and $\beta=0.9$. Means $\lambda_i^t$ of Poisson demand in each period $t=1,\ldots,T$ were integer numbers uniformly distributed between 5 and 10. 

We fixed $\alpha=0.9$, $\vartheta=0.05$ and $M=10$; recall that $M$ is the number of independently generated instances of $\mathcal{P}_{\text{lb}}$ and $\mathcal{P}_{\text{ub}}$ used for computing profit/cost upper and lower bounds as illustrated in Section \ref{sec:scop}. This led to a sample size of 209 for the SSMKP and of 370 for the SSLSP (Definition \ref{sample_size_alpha_theta_solution_set_approx}).

Due to the small size of the SSMKP instances, we managed to obtain optimal solutions by exhaustive enumeration, i.e. we generated all possible assignment and then checked feasibility and expected total profit of each of them via Monte Carlo simulation. The number of Monte Carlo runs was set to guarantee a margin of error of $\vartheta/10$ with a confidence level of $0.9$, in such a way as to ensure an approximation error negligible with respect to the chosen $\vartheta$. SSLSP instances can be solved to optimality by using a deterministic equivalent mixed integer linear programming model \cite{citeulike:7292595}. In our analysis, we can therefore compare results obtained with our approach against the true optimal solutions.

In Fig. \ref{fig:alpha_theta_opt_solutions}, for each instance, we plotted upper and lower bound obtained for its optimal profit (SSMKP) or cost (SSLSP). For clarity, the interval has been normalised by using the profit/cost of the true optimal solution as a normalisation factor, so that value one in the graph denotes the true optimal profit/cost. The confidence level achieved by using our approach is generally higher than the prescribed $\alpha$. In fact, despite $\alpha$ being set to 0.9, over the hundred instances analysed, the cost confidence interval did not cover the true optimal cost only in one case (SSMKP, instance 21). This is due to the conservative nature of our approach, as already discussed in Section \ref{sec:feasibility}. 

We believe the fluctuations in the size of optimality gaps observed in Fig. \ref{fig:alpha_theta_opt_solutions} for the SSMKP may be related to the fact that this problem features 0-1 integer variables. Depending on the specific instance being solved, different sets of samples may lead to assignments in which ``high value'' objects belonging to the true optimal solution of the problem are not selected. This may lead to larger optimality gaps than those observed for other instances in which the optimal solution is less sensitive to random fluctuations produced by the sampling process. 
 
\begin{figure}[ht]
\begin{center}
\[\begin{array}{cc}
\resizebox{0.48\columnwidth}{!}{ 
\begin{tikzpicture}
\begin{axis}[
title={SSMKP},
xtick = {1,10,20,30,40,50},
xlabel=Instance,
ylabel=Normalised profit,
yticklabel style={
        /pgf/number format/fixed,
        /pgf/number format/precision=5
},
scaled y ticks=false
]

\addplot[smooth,densely dotted,mark=none,
domain=-1:51,samples=40]
{1};

\addplot+[black,
mark=point,
only marks,
error bars/.cd,
y dir=minus,y explicit]
coordinates {
(	1	,	1	)	+-	(0,	0.009689376	)
(	2	,	1	)	+-	(0,	0.121389863	)
(	3	,	1	)	+-	(0,	0.008010189	)
(	4	,	1	)	+-	(0,	0.011535649	)
(	5	,	1	)	+-	(0,	0.010890018	)
(	6	,	1	)	+-	(0,	0.004361218	)
(	7	,	1	)	+-	(0,	0.009035399	)
(	8	,	1	)	+-	(0,	0.110444576	)
(	9	,	1	)	+-	(0,	0.01243542	)
(	10	,	1	)	+-	(0,	0.020960773	)
(	11	,	1	)	+-	(0,	0.009902942	)
(	12	,	1	)	+-	(0,	0.012867987	)
(	13	,	1	)	+-	(0,	0.014991301	)
(	14	,	1	)	+-	(0,	0.150516871	)
(	15	,	1	)	+-	(0,	0.021729527	)
(	16	,	1	)	+-	(0,	0.006981601	)
(	17	,	1	)	+-	(0,	0.009308936	)
(	18	,	1	)	+-	(0,	0.004570477	)
(	19	,	1	)	+-	(0,	0.013023245	)
(	20	,	1	)	+-	(0,	0.125981825	)
(	21	,	1	)	+-	(0,	-0.000187619	)
(	22	,	1	)	+-	(0,	0.00644165	)
(	23	,	1	)	+-	(0,	0.170080344	)
(	24	,	1	)	+-	(0,	0.012527795	)
(	25	,	1	)	+-	(0,	0.014991542	)
(	26	,	1	)	+-	(0,	0.010177072	)
(	27	,	1	)	+-	(0,	0.008551275	)
(	28	,	1	)	+-	(0,	0.005970499	)
(	29	,	1	)	+-	(0,	0.109675686	)
(	30	,	1	)	+-	(0,	0.113178055	)
(	31	,	1	)	+-	(0,	0.01748577	)
(	32	,	1	)	+-	(0,	0.011061922	)
(	33	,	1	)	+-	(0,	0.014951763	)
(	34	,	1	)	+-	(0,	0.0173359	)
(	35	,	1	)	+-	(0,	0.011209656	)
(	36	,	1	)	+-	(0,	0.011043827	)
(	37	,	1	)	+-	(0,	0.020124338	)
(	38	,	1	)	+-	(0,	0.009764761	)
(	39	,	1	)	+-	(0,	0.008186626	)
(	40	,	1	)	+-	(0,	0.012937778	)
(	41	,	1	)	+-	(0,	0.009374519	)
(	42	,	1	)	+-	(0,	0.00535904	)
(	43	,	1	)	+-	(0,	0.008918388	)
(	44	,	1	)	+-	(0,	0.005411962	)
(	45	,	1	)	+-	(0,	0.014357858	)
(	46	,	1	)	+-	(0,	0.019304972	)
(	47	,	1	)	+-	(0,	0.120707623	)
(	48	,	1	)	+-	(0,	0.008847918	)
(	49	,	1	)	+-	(0,	0.005800859	)
(	50	,	1	)	+-	(0,	0.007979165	)
};

\addplot+[black,
mark=point,
only marks,
error bars/.cd,
y dir=plus,y explicit]
coordinates {
(	1	,	1	)	+-	(0,	0.139633438	)
(	2	,	1	)	+-	(0,	0.010601886	)
(	3	,	1	)	+-	(0,	0.016636474	)
(	4	,	1	)	+-	(0,	0.011878386	)
(	5	,	1	)	+-	(0,	0.01033615	)
(	6	,	1	)	+-	(0,	0.021896207	)
(	7	,	1	)	+-	(0,	0.011527367	)
(	8	,	1	)	+-	(0,	0.007283032	)
(	9	,	1	)	+-	(0,	0.017155044	)
(	10	,	1	)	+-	(0,	0.007311475	)
(	11	,	1	)	+-	(0,	0.004681728	)
(	12	,	1	)	+-	(0,	0.015319601	)
(	13	,	1	)	+-	(0,	0.01515816	)
(	14	,	1	)	+-	(0,	0.01517206	)
(	15	,	1	)	+-	(0,	0.014011408	)
(	16	,	1	)	+-	(0,	0.029307619	)
(	17	,	1	)	+-	(0,	0.018659067	)
(	18	,	1	)	+-	(0,	0.018422925	)
(	19	,	1	)	+-	(0,	0.013515367	)
(	20	,	1	)	+-	(0,	0.016345394	)
(	21	,	1	)	+-	(0,	0.153001486	)
(	22	,	1	)	+-	(0,	0.024557207	)
(	23	,	1	)	+-	(0,	0.02168767	)
(	24	,	1	)	+-	(0,	0.022286644	)
(	25	,	1	)	+-	(0,	0.015092042	)
(	26	,	1	)	+-	(0,	0.137466184	)
(	27	,	1	)	+-	(0,	0.013464176	)
(	28	,	1	)	+-	(0,	0.017811662	)
(	29	,	1	)	+-	(0,	0.010511252	)
(	30	,	1	)	+-	(0,	0.014695067	)
(	31	,	1	)	+-	(0,	0.009230531	)
(	32	,	1	)	+-	(0,	0.017602347	)
(	33	,	1	)	+-	(0,	0.019544713	)
(	34	,	1	)	+-	(0,	0.008540532	)
(	35	,	1	)	+-	(0,	0.023375119	)
(	36	,	1	)	+-	(0,	0.011756669	)
(	37	,	1	)	+-	(0,	0.011937739	)
(	38	,	1	)	+-	(0,	0.01927978	)
(	39	,	1	)	+-	(0,	0.015963658	)
(	40	,	1	)	+-	(0,	0.016935311	)
(	41	,	1	)	+-	(0,	0.014777846	)
(	42	,	1	)	+-	(0,	0.188833348	)
(	43	,	1	)	+-	(0,	0.018377036	)
(	44	,	1	)	+-	(0,	0.15632719	)
(	45	,	1	)	+-	(0,	0.01742092	)
(	46	,	1	)	+-	(0,	0.016455438	)
(	47	,	1	)	+-	(0,	0.01265442	)
(	48	,	1	)	+-	(0,	0.014429074	)
(	49	,	1	)	+-	(0,	0.020392372	)
(	50	,	1	)	+-	(0,	0.177809475	)
};
\end{axis}
\end{tikzpicture}
}
&
\resizebox{0.48\columnwidth}{!}{ 
\begin{tikzpicture}
\begin{axis}[
title={SSLSP},
xtick = {1,10,20,30,40,50},
xlabel=Instance,
ylabel=Normalised cost,
yticklabel style={
        /pgf/number format/fixed,
        /pgf/number format/precision=5
},
scaled y ticks=false
]

\addplot[smooth,densely dotted,mark=none,
domain=-1:51,samples=40]
{1};

\addplot+[black,
mark=circle,
only marks,
error bars/.cd,
y dir=minus,y explicit]
coordinates {
(	1	,	1	)	+-	(0,	0.090909091	)
(	2	,	1	)	+-	(0,	0.118421053	)
(	3	,	1	)	+-	(0,	0.111111111	)
(	4	,	1	)	+-	(0,	0.093333333	)
(	5	,	1	)	+-	(0,	0.115384615	)
(	6	,	1	)	+-	(0,	0.133333333	)
(	7	,	1	)	+-	(0,	0.1	)
(	8	,	1	)	+-	(0,	0.13253012	)
(	9	,	1	)	+-	(0,	0.113924051	)
(	10	,	1	)	+-	(0,	0.106666667	)
(	11	,	1	)	+-	(0,	0.126760563	)
(	12	,	1	)	+-	(0,	0.103896104	)
(	13	,	1	)	+-	(0,	0.115384615	)
(	14	,	1	)	+-	(0,	0.121621622	)
(	15	,	1	)	+-	(0,	0.119047619	)
(	16	,	1	)	+-	(0,	0.108108108	)
(	17	,	1	)	+-	(0,	0.102941176	)
(	18	,	1	)	+-	(0,	0.105263158	)
(	19	,	1	)	+-	(0,	0.097222222	)
(	20	,	1	)	+-	(0,	0.109589041	)
(	21	,	1	)	+-	(0,	0.105263158	)
(	22	,	1	)	+-	(0,	0.131578947	)
(	23	,	1	)	+-	(0,	0.118421053	)
(	24	,	1	)	+-	(0,	0.115384615	)
(	25	,	1	)	+-	(0,	0.112676056	)
(	26	,	1	)	+-	(0,	0.12	)
(	27	,	1	)	+-	(0,	0.115384615	)
(	28	,	1	)	+-	(0,	0.113924051	)
(	29	,	1	)	+-	(0,	0.088235294	)
(	30	,	1	)	+-	(0,	0.093333333	)
(	31	,	1	)	+-	(0,	0.116883117	)
(	32	,	1	)	+-	(0,	0.08974359	)
(	33	,	1	)	+-	(0,	0.123287671	)
(	34	,	1	)	+-	(0,	0.098591549	)
(	35	,	1	)	+-	(0,	0.118421053	)
(	36	,	1	)	+-	(0,	0.106666667	)
(	37	,	1	)	+-	(0,	0.12	)
(	38	,	1	)	+-	(0,	0.134146341	)
(	39	,	1	)	+-	(0,	0.109589041	)
(	40	,	1	)	+-	(0,	0.12987013	)
(	41	,	1	)	+-	(0,	0.103896104	)
(	42	,	1	)	+-	(0,	0.094594595	)
(	43	,	1	)	+-	(0,	0.103896104	)
(	44	,	1	)	+-	(0,	0.116883117	)
(	45	,	1	)	+-	(0,	0.106666667	)
(	46	,	1	)	+-	(0,	0.101265823	)
(	47	,	1	)	+-	(0,	0.106666667	)
(	48	,	1	)	+-	(0,	0.135802469	)
(	49	,	1	)	+-	(0,	0.12	)
(	50	,	1	)	+-	(0,	0.115942029	)
};

\addplot+[black,
mark=circle,
only marks,
error bars/.cd,
y dir=plus,y explicit]
coordinates {
(	1	,	1	)	+-	(0,	0.142857143	)
(	2	,	1	)	+-	(0,	0.171052632	)
(	3	,	1	)	+-	(0,	0.166666667	)
(	4	,	1	)	+-	(0,	0.173333333	)
(	5	,	1	)	+-	(0,	0.128205128	)
(	6	,	1	)	+-	(0,	0.146666667	)
(	7	,	1	)	+-	(0,	0.185714286	)
(	8	,	1	)	+-	(0,	0.156626506	)
(	9	,	1	)	+-	(0,	0.164556962	)
(	10	,	1	)	+-	(0,	0.133333333	)
(	11	,	1	)	+-	(0,	0.14084507	)
(	12	,	1	)	+-	(0,	0.155844156	)
(	13	,	1	)	+-	(0,	0.166666667	)
(	14	,	1	)	+-	(0,	0.162162162	)
(	15	,	1	)	+-	(0,	0.119047619	)
(	16	,	1	)	+-	(0,	0.148648649	)
(	17	,	1	)	+-	(0,	0.132352941	)
(	18	,	1	)	+-	(0,	0.144736842	)
(	19	,	1	)	+-	(0,	0.152777778	)
(	20	,	1	)	+-	(0,	0.164383562	)
(	21	,	1	)	+-	(0,	0.171052632	)
(	22	,	1	)	+-	(0,	0.105263158	)
(	23	,	1	)	+-	(0,	0.184210526	)
(	24	,	1	)	+-	(0,	0.166666667	)
(	25	,	1	)	+-	(0,	0.154929577	)
(	26	,	1	)	+-	(0,	0.16	)
(	27	,	1	)	+-	(0,	0.166666667	)
(	28	,	1	)	+-	(0,	0.139240506	)
(	29	,	1	)	+-	(0,	0.132352941	)
(	30	,	1	)	+-	(0,	0.133333333	)
(	31	,	1	)	+-	(0,	0.168831169	)
(	32	,	1	)	+-	(0,	0.153846154	)
(	33	,	1	)	+-	(0,	0.164383562	)
(	34	,	1	)	+-	(0,	0.183098592	)
(	35	,	1	)	+-	(0,	0.144736842	)
(	36	,	1	)	+-	(0,	0.16	)
(	37	,	1	)	+-	(0,	0.133333333	)
(	38	,	1	)	+-	(0,	0.170731707	)
(	39	,	1	)	+-	(0,	0.191780822	)
(	40	,	1	)	+-	(0,	0.181818182	)
(	41	,	1	)	+-	(0,	0.181818182	)
(	42	,	1	)	+-	(0,	0.175675676	)
(	43	,	1	)	+-	(0,	0.155844156	)
(	44	,	1	)	+-	(0,	0.155844156	)
(	45	,	1	)	+-	(0,	0.133333333	)
(	46	,	1	)	+-	(0,	0.139240506	)
(	47	,	1	)	+-	(0,	0.173333333	)
(	48	,	1	)	+-	(0,	0.148148148	)
(	49	,	1	)	+-	(0,	0.146666667	)
(	50	,	1	)	+-	(0,	0.130434783	)
};
\end{axis}
\end{tikzpicture}
}
\end{array}\]
\end{center}
\caption{Normalised profit/cost upper and lower bounds for fifty SSMKP and SSLSP instances; a value of 1 denotes the true optimal profit/cost.}
\label{fig:alpha_theta_opt_solutions}
\end{figure}
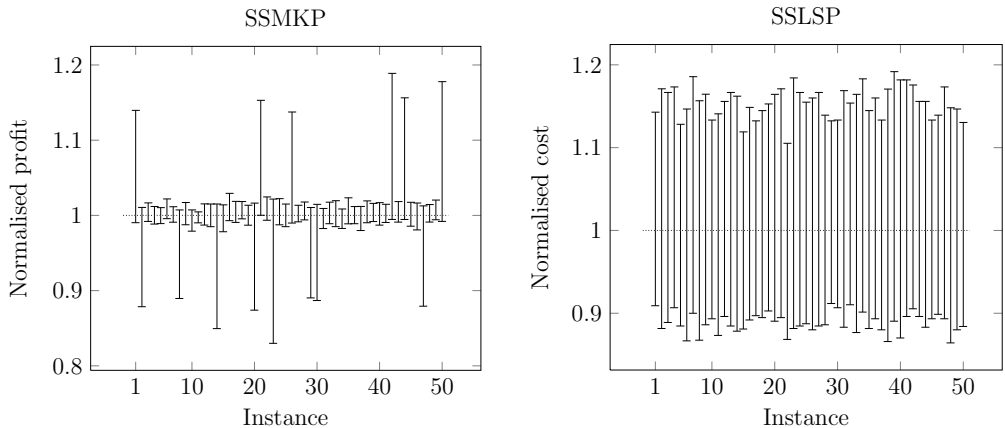

\begin{figure}[h!]
\begin{center}
\[\begin{array}{cc}
\resizebox{0.48\columnwidth}{!}{ 
\begin{tikzpicture}
\begin{axis}[
	legend style={at={(0.5,-0.2)},
	anchor=north,legend columns=-1},
xtick = {0.01,0.025,0.05},
xlabel=$\vartheta$,
ylabel=Optimality gap \%,
yticklabel style={
        /pgf/number format/fixed,
        /pgf/number format/precision=5
},
scaled y ticks=false,
xticklabel style={
        /pgf/number format/fixed,
        /pgf/number format/precision=5
},
scaled x ticks=false
]
\addplot[black]
	coordinates {(0.05,6.234535677) (0.025,2.84916914) (0.01,0.668696219)};
\addplot[black,dashed]
	coordinates {(0.05,26.71110974) (0.025,12.60315256) (0.01,5.013945069)};
\legend{SSMKP,SSLSP}
\end{axis}
\end{tikzpicture}
}
&
\resizebox{0.48\columnwidth}{!}{ 
\begin{tikzpicture}
\begin{axis}[
	legend style={at={(0.5,-0.2)},
	anchor=north,legend columns=-1},
xtick = {0.1,0.2,0.3},
xlabel=$\vartheta$,
ylabel=Optimality gap (periods),
yticklabel style={
        /pgf/number format/fixed,
        /pgf/number format/precision=5
},
scaled y ticks=false,
xticklabel style={
        /pgf/number format/fixed,
        /pgf/number format/precision=5
},
scaled x ticks=false
]
\addplot[black]
	coordinates {(0.1,0.84) (0.2,1.42) (0.3,2.06)};
\legend{SMPSP}
\end{axis}
\end{tikzpicture}
}
\end{array}\]
\end{center}
\caption{Average optimality gap for different values of $\vartheta$}
\label{fig:opt_gap_vartheta}
\end{figure}
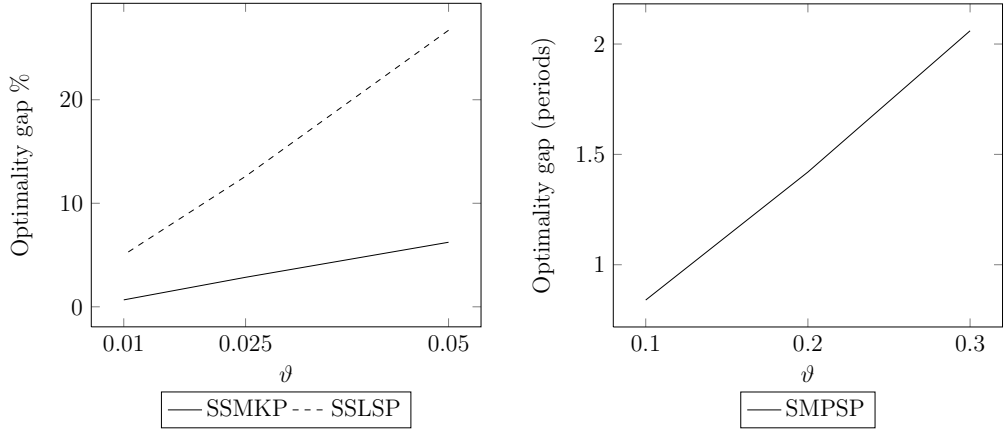

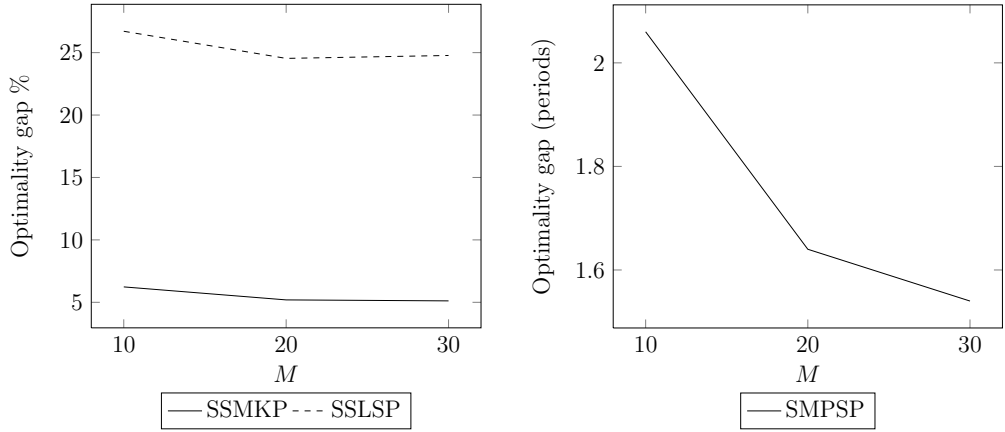
\begin{figure}[h!]
\begin{center}
\[\begin{array}{cc}
\resizebox{0.48\columnwidth}{!}{ 
\begin{tikzpicture}
\begin{axis}[
	legend style={at={(0.5,-0.2)},
	anchor=north,legend columns=-1},
xtick = {10,20,30},
xlabel=$M$,
ylabel=Optimality gap \%,
yticklabel style={
        /pgf/number format/fixed,
        /pgf/number format/precision=5
},
scaled y ticks=false,
xticklabel style={
        /pgf/number format/fixed,
        /pgf/number format/precision=5
},
scaled x ticks=false
]
\addplot[black]
	coordinates {(10,6.234535677) (20,5.193449725) (30,5.114829553)};
\addplot[black,dashed]
	coordinates {(10,26.71110974) (20,24.53841982) (30,24.77254963)};
\legend{SSMKP,SSLSP}
\end{axis}
\end{tikzpicture}
}
&
\resizebox{0.48\columnwidth}{!}{ 
\begin{tikzpicture}
\begin{axis}[
	legend style={at={(0.5,-0.2)},
	anchor=north,legend columns=-1},
xtick = {10,20,30},
xlabel=$M$,
ylabel=Optimality gap (periods),
yticklabel style={
        /pgf/number format/fixed,
        /pgf/number format/precision=5
},
scaled y ticks=false,
xticklabel style={
        /pgf/number format/fixed,
        /pgf/number format/precision=5
},
scaled x ticks=false
]
\addplot[black]
	coordinates {(10,2.06) (20,1.64) (30,1.54)};
\legend{SMPSP}
\end{axis}
\end{tikzpicture}
}
\end{array}\]
\end{center}
\caption{Average optimality gap for different values of $M$}
\label{fig:opt_gap_M}
\end{figure}

Finally, we included in the analysis randomly generated instances of the SMPSP formulated as an SCOP in which the objective is to minimise the latest start time. In these instances $K=5$ and $\beta=0.6$; $r_k$ and $d_k$, which represent job $k$ release time and deadline, were all set to $0$ and $20$, respectively. Capacity requirements $c_k$ where generated as integer numbers uniformly distributed between 1 and 3. Expected task durations $\lambda_k$ were generated as uniformly distributed numbers between 1 and 5; the maximum number of processors $P$ was set to 5. 

In Fig. \ref{fig:opt_gap_vartheta} we analysed the behaviour of the optimality gap when $\alpha=0.9$, $M=10$ and $\vartheta$ varies. The sample size ranges as follows: from 209 ($\vartheta=0.05$) to 5838 ($\vartheta=0.01$) for the SSMKP; from 370 ($\vartheta=0.05$) to 6350 ($\vartheta=0.01$) for the SSLSP; and from 14 ($\vartheta=0.3$) to 114 ($\vartheta=0.1$) for the SMPSP. For each value of $\vartheta$ considered, we solved 50 different instances of the SSMKP, SSLSP and SMPSP, and we computed the average optimality gap over this pool of instances. The average optimality gap for the SSMKP and the SSLSP is reported in percentage of the true optimal solution. For the case of the SMPSP unfortunately we were not able to compute the true optimal plan, therefore we reported the optimality gap in absolute terms; since we are minimising the latest start time, we expressed the optimality gap in expected number of periods. Note that since $\alpha$ and $\vartheta$ are linked to the number of samples generated by the relation in Definition \ref{sample_size_alpha_theta_solution_set_approx}, similar plots may be obtained by varying $\alpha$ and keeping $\vartheta$ fixed. 

In Fig. \ref{fig:opt_gap_M} we carried out a similar analysis by keeping $\vartheta$ fixed to 0.05 (SSMKP, SSLSP) and to 0.3 (SSMKP) and by varying $M$. 

\subsection{Computational efficiency}\label{sec:comp_efficiency}

In this section we reflect on the computational complexity and on the scalability of our approach. 

\subsubsection{Computational complexity}

The computational complexity of SCSPs has been discussed in several works \cite{DBLP:conf/ecai/Walsh02,tmw2006,citeulike:10689458}; in particular we direct the interested reader to \cite{citeulike:13582209, citeulike:11058108}, which provide comprehensive overviews on the complexity of stochastic programs. 

Multi-stage stochastic programming with discretely distributed decision-dependent random variables is PSPACE-hard \cite{citeulike:11058108}; the result follows from the PSPACE-hardness of the problem ``decision-making under uncertainty'' in \cite{citeulike:13582254}. SCSPs are PSPACE-complete in general if random variables are defined on discrete supports \cite{tmw2006}. However, as pointed out in \cite{citeulike:11058108}, the complexity of the ``standard'' multi-stage stochastic programming problem, in which distributions are independent of decisions taken in earlier stages, remains open; the authors in \cite{citeulike:11058108} conjecture that this is also PSPACE-hard.

The advantage of sampled SCSPs over generic SCSPs is that sampled SCSPs always comprise a finite number of scenarios whose number is determined by Definition \ref{sample_size_alpha_theta_solution_set_approx}, which establishes a relationship among $\alpha$, $\vartheta$, and $N$. A decision maker is then free to fix a pair of these values and to derive the remaining one. In principle, one may fix a priori the number of samples $N$ --- rather than the confidence level $\alpha$ or the error threshold $\vartheta$ --- and sacrifice precision for efficiency. This will not make the sampled SCSP fixed-parameter tractable in general --- in \cite{citeulike:10689458} the authors proved that maintaining GAC on a global chance constraint can be intractable even when maintaining GAC on the corresponding deterministic version of that constraint is tractable --- but it may reduce its complexity from PSPACE to NP-hard. 

\subsubsection{Scalability}

To illustrate the scalability of our approach with respect to other state-of-the-art approaches to SCSPs we employ, once more, the SSMKP. Note that this problem is similar to the one discussed in \citep[][Section 8.3]{citeulike:10689458}. In Section 9.4 of the same work, it was discussed that --- for this class of problems --- even when profits and weights of the objects are defined on a support that comprises only two values, a scenario-based formulation would end up comprising $2^{20}$ scenarios and a solver such as Choco would run out of memory. It was then shown that the complete approach discussed in that work could solve an instance comprising 10 objects in about an hour on average. 

We fixed $\alpha=0.9$, $\vartheta=0.01$ and $M=10$. We solved ten instances of the SSMKP randomly generated as discussed in Section \ref{sec:optimality}, but now comprising $N=20$ rather than ten objects; this led to a sample size of 6916. In Table \ref{tab:SSMKP} we report the optimality gap and the runtime for each of these instances as well as the runtime in hours.

\noindent
\begin{table}
\centering
\begin{tabular}{lrr}
Instance&Optimality gap \%&Runtime (hours)\\
\hline
1	&	0.49	&	1.65\\
2	&	0.66	&	0.98\\
3	&	0.57	&	1.20\\
4	&	0.71	&	0.65\\
5	&	0.53	&	1.77\\
6	&	0.51	&	2.18\\
7	&	0.50	&	2.22\\
8	&	0.55	&	0.77\\
9	&	0.44	&	2.63\\
10	&	0.70	&	0.80\\
\hline					
Mean	&	0.57	&	1.49	
\end{tabular}
\caption{SSMKP. Larger instances comprising $N=20$ objects}
\label{tab:SSMKP}
\end{table}
    
Finally, we fixed $\alpha=0.9$, $\vartheta=0.1$ and $M=10$, and we solved larger instances of the SMPSP formulated as an SCOP. These instances comprise ten jobs (i.e. $K=10$). $r_k$ and $d_k$ were now set to $0$ and $30$, respectively; this led to a sample size of 142. In Table \ref{tab:SMPSP} we report upper and lower bound for the latest start time associated with each of these instances, as well as and the runtime in hours. 

\begin{table}
\centering
\begin{tabular}{lrrr}
&\multicolumn{2}{c}{Latest Start Time}\\
Instance&LB&UB&Runtime (hours)\\
\hline
1	&	13	&	15	&	1.50	\\
2	&	14	&	16	&	0.68	\\
3	&	11	&	12	&	0.72	\\
4	&	13	&	14	&	3.20	\\
5	&	17	&	21	&	0.31	\\
6	&	17	&	19	&	3.27	\\
7	&	18	&	20	&	9.08	\\
8	&	13	&	14	&	0.43	\\
9	&	13	&	14	&	0.16	\\
10	&	17	&	20	&	11.11\\
\hline							
Mean	&	&		&	3.05	
\end{tabular}
\caption{SMPSP. Larger instances comprising $K=10$ jobs}
\label{tab:SMPSP}
\end{table}

It is clear that it would be impossible to directly use the approach in \citep{citeulike:10689458} to model these SSMKP instances, as random variables follow a Poisson distribution and therefore have infinite values in their support. Even if one discretises these supports, e.g. by reducing them to only two values, the resulting SCSP would feature millions of scenarios.

However, one may argue that, in the case of the SSMKP, we are analysing a problem that could be analysed by brute force. In other words, one may as well generate all $2^{20}$ possible assignments for the decision variables and then analytically check the feasibility of each. Unfortunately, it is clear that this is not possible for the two-stage SMPSP just analysed, which features a much larger search space.

These results therefore demonstrate that the discussion in this work provides a viable means for scaling up the approach in \citep{citeulike:10689458}. 

\section{Related works}\label{sec:rel_works}

Confidence-based optimisation was originally introduced in \cite{citeulike:13269658}. In this work, the authors discuss an application of this methodology in the context of a well-known stochastic inventory control problem. Our work extends the discussion presented there to generic SCSPs and SCOPs by introducing a more general notion of confidence-based reasoning based on two novel concepts: ($\alpha,\vartheta$)-solutions and ($\alpha,\vartheta$)-solution sets. In the context of stochastic modeling and optimisation, as discussed, these tools can be employed to find approximate solutions that possess given statistical properties. 
 
\subsection{Related works in stochastic programming}

In operations research, and particularly in stochastic programming, the state-of-the-art technique that applies sampling in combinatorial optimisation is the sample average approximation (SAA) method \cite{ksh01}. This is a Monte Carlo simulation-based approach to stochastic discrete optimisation problem. This method replaces the actual distribution of random variables in the combinatorial problem of interest by an empirical distribution obtained via sampling. The obtained ``sample average optimisation problem'' is then solved and the procedure is repeated multiple times until a given termination criterion is satisfied.
The authors in \cite{ksh01} focus on stochastic programs with expected value objectives and discuss convergence rate and stopping rules. In \cite{Ahmed02thesample} the authors extend their analysis to two-stage stochastic programs with integer recourse; for this latter class of problems \cite{citeulike:541752} carry out a post-hoc computationally intensive analysis of the quality of solutions obtained via SAA. Extensions to problems with expected value constraints, e.g. conditional value-at-risk constraints, were discussed in \cite{citeulike:13321118}. However, none of these works investigated the case in which the problem of interest include chance constraints. As \cite{citeulike:10237336} remarks, there are formulations of stochastic programming problems that incorporate expectations of penalised constraints in the objective function as a penalty terms. These problems can be solved efficiently since they simply require continuous variables for modelling penalties and they do not require any additional binary variable. However, this modelling approach does not address the issue of finding or approximating feasible or optimal solutions to a chance constrained problem \cite[][p. 950]{citeulike:13324019}. 

SAA methods for problems comprising a single chance constraint were discussed in \cite{ahm_shap_08,citeulike:9453433,citeulike:4176929}. In \cite{ahm_shap_08} the authors summarise convergence properties (Section 2.1) and post-hoc solution validation strategies (Section 2.2). They remark that ``based on this [convergence] analysis, we can compute a priori the sample size required in the SAA problem so that it produces a feasible solution to the true problem with high probability (typically such estimates of a required sample size are quite conservative).'' The convergence analysis the authors refer to was originally conducted in \cite{citeulike:9453433} and it shows asymptotical convergence properties based on inequalities such as Chernoff's \cite{citeulike:9389607} or Hoeffding's \cite{citeulike:3392582}, which are known to be conservative bounds. Their analysis is conducted under the assumption that the feasible region is finite, since the sample size determined via the aforementioned convergence properties grows linearly in the size of the feasible region \cite[][p. 683]{citeulike:9453433}. Extensions of the analysis in \cite{citeulike:9453433} to the case of multiple chance constraints were illustrated in \cite{citeulike:10237336}. This latter work is similar to those just discussed, since once more the analysis is based on the above inequalities and the sample size depends on the size of the feasible region. More recently, \cite{citeulike:13324051} investigated the relations between chance constrained and penalty function problems under discrete distributions. This analysis extended a number of previous works that analysed this relation under continuous distribution. However, the authors explicitly remark that ``our goal is not to show that the penalty problems are able to generate optimal values and solutions of chance constrained problems.'' Instead they compare the problems with focus on asymptotic equivalence of optimal values and corresponding convergence of optimal solutions.

After surveying the existing literature on SAA, the first important remark is that in none of the above works can we find concepts that resemble those of ($\alpha$,$\vartheta$)-solution and ($\alpha$,$\vartheta$)-solution set, which are unique to confidence-based reasoning \cite{citeulike:13269658}. This is a subtle conceptual difference that should not be overlooked. The aim of SAA is to find an assignment that, with prescribed confidence probability $\alpha$, is a solution to the original problem; see e.g. \cite[][Section 3.1]{citeulike:4176929}. In other words, in SAA the decision maker does not fix any a priori tolerated estimation error $\vartheta$. To ensure that the solution of the sampled problem is feasible with respect to the original problem with sufficiently high probability, in SAA the threshold $\beta$ associated with chance constraints in the sampled problem is increased by a factor $\vartheta$, which however is not explicitly interpreted as an error tolerance threshold in a statistical sense, although in practice it is used as such. In \cite[][p. 407]{citeulike:4176929}, the authors point out that, for a fixed $\alpha$ and for a given threshold $\beta$, ``it is not clear what the best choices for the sample size and $\vartheta$ are,'' since they believe this is a problem-dependent issue that should be addressed numerically. This statement demonstrates the aforementioned fundamental difference. By introducing the two concepts of ($\alpha$,$\vartheta$)-solution and ($\alpha$,$\vartheta$)-solution set, we suggest that a decision maker may --- in line with established practices in statistics --- interpret $\vartheta$ as an error tolerance threshold and fix a priori, together with the confidence level $\alpha$, either the sample size (on the basis of the available observations) or $\vartheta$ (on the basis of the estimation error that can be tolerated); finally, the parameter that has not been fixed should be derived via the analysis we presented. In summary, the difference lies in the interpretation. Confidence-based reasoning aims to find a solution that, with confidence $\alpha$, satisfies the chance constraints in the original problem within the given error tolerance $\vartheta$. In addition to this important semantic difference, we should mention that our analysis is based on the exact Clopper-Pearson confidence interval, and not on conservative bounds such as Chernoff's or Hoeffding's inequalities. Finally, our approximation strategy for ($\alpha$,$\vartheta$)-solution sets leads to a sample size 
(Definition \ref{sample_size_alpha_theta_solution_set_approx}) that is independent of the number of assignments in the feasible region; a major difference from all other methods surveyed so far.  

To contrast our approach with respect to other existing state-of-the-art approaches, one may consider the stochastic vehicle routing problem with time windows discussed in \cite[][Section 4]{citeulike:10237336}. It is possible to apply our analysis to the instances discussed in \cite[][Table 1]{citeulike:10237336}, by converting the parameters used in SAA and setting the confidence level $\alpha=0.99$, the error threshold $\vartheta=0.05$ and the chance constraint thresholds $\beta=0.95$. For the instances with 10 customer orders, the sample size prescribed by our approach is 429 for the model with a single chance constraint and 490 for the model with three chance constraints. For the instances with 50 customer orders, the sample size prescribed by our approach is 608 for the model with a single chance constraint and 669 for the model with three chance constraints. Not only are these sample sizes orders of magnitude smaller than the ones suggested in \cite{citeulike:10237336}, which range from 200 thousand up to 32 million; but most importantly they do not depend on the number of vessels used or the size of the time windows; in fact, according to Definition \ref{sample_size_alpha_theta_solution_set_approx}, they only depend on the number of random variables and chance constraints in the model. Of course, as the authors in \cite{citeulike:10237336} remark, finding an exact solution to a scenario-based model with 32 million scenarios is unrealistic. For this reason, they suggest to adopt heuristic solution methods, e.g. tabu search. As demonstrated in our computational study, our approximate ($\alpha$,$\vartheta$)-solution sets represent a viable alternative to the use of heuristics on instances featuring very large sample sizes. 

\subsection{Related works in constraint programming}

A detailed discussion on hybrid CP/AI/OR approaches for decision making under uncertainty can be found in \cite{BrownMiguel06, hnich_et_al_survey_2010}. We direct the reader to these two references for further details on existing works in this research area. We next briefly survey key relevant references.
Efforts that try to extend classical CSP framework to
incorporate uncertainty have been influenced by works that 
originated in different fields, namely {\em chance-constrained
programming} \cite{cc63} and {\em stochastic programming}
\cite{bl97}. 
To the best of our knowledge the first work that tries to create a bridge 
between Stochastic Programming and Constraint 
Programming is by Benoist et al. \cite{bbcr01}.
Search and consistency strategies, namely a backtracking
algorithm, a forward checking procedure \cite{DBLP:conf/ecai/Walsh02} and
an arc-consistency \cite{DBLP:conf/cp/BalafoutisS06} algorithm 
have been proposed for SCSPs.
A scenario-based approach for building up constraint programming 
models of SCSPs was proposed by Tarim et al. \cite{tmw2006}. 
In the same work a fully featured language --- Stochastic OPL --- for 
modeling SCSPs was also proposed. 
In \cite{1244077} the authors introduce new algorithms for solving 
multi-objective stochastic problems are proposed.
Global chance constraints were introduced first in \cite{ros_cons_07},
and bring together the reasoning power of global constraints from CP
and the expressive power of chance constraints from SP. A general
purpose approach for filtering global chance constraints is proposed
in \cite{DBLP:conf/cp/HnichRTP09,citeulike:10689458}. This approach is able to reuse
existing propagators available for the respective deterministic global
constraint which corresponds to a given global chance constraint when
all the random variables are replaced by constant parameters. 
In \cite{DBLP:conf/cp/RossiTHP08} the authors discuss some
possible strategies to perform cost-based filtering for certain
classes of SCOPs. These strategies exploit well-known
inequalities borrowed from SP and used to compute valid bounds for any
given SCOP that respects some mild assumptions. 
Unfortunately, above approaches operate under the assumption 
that the number of scenarios must be finite, 
otherwise a solution cannot be expressed as a finite
number of possible decisions. Furthermore,
these approaches do not scale well. Even problems having a limited number of stochastic
variables with large support immediately produce policy trees whose size makes impractical
the use of a complete method. 
In \cite{tmw2006} the authors employed sampling in order to reduce the
number of scenarios considered for a given stochastic constraint program 
and produce a solution in reasonable time. Nevertheless, this approach
does not provide any optimality/feasibility guarantee for the solution produced.
Heuristic approaches such as the one in \cite{DBLP:conf/cp/PrestwichTRH09}, 
in which a neural network is employed in order to encode a policy function, 
suffer from the same limitation and from lack of modularity. Stochastic sampling in the context
of Stochastic Boolean Satisfiability was discussed in \cite{DBLP:journals/jar/LittmanMP01};
forward sampling \cite{bbcr01}  and 
sample aggregation \cite{DBLP:conf/cpaior/HentenryckBV06} 
are two other techniques that have been employed to solve SCSPs.
Nevertheless, none of these approaches introduce a concept that
resembles that of ($\alpha$,$\vartheta$)-solution. 
Probably, the work discussed in \cite{DBLP:conf/icalp/KatrielKU07} 
represents the closest attempt to provide some sort of guarantees
for a stochastic constraint satisfaction problem. Nevertheless, this
work is focused on a specific problem --- a two-stage stochastic 
matching problem --- and it does not propose a generic 
approach for solving SCSPs. Finally, another closely related work is
\cite{DBLP:journals/jair/BeckW07}, which discusses sample-based approaches
to job shop scheduling with probabilistic durations; however, like in the previous case,
the approach proposed is focused on a specific problem and not on solving 
generic SCSPs.

\section{Conclusions}\label{sec:con}

We proposed a framework for exploiting sampling in order to solve 
SCSPs that include random variables over a continuous or very large discrete support. 
Our framework is based on a number of novel concepts: sampled SCSPs, $(\alpha,\vartheta)$-solutions and $(\alpha,\vartheta)$-solution sets. We employed statistical estimation to determine if a given assignment 
is consistent with respect to a given set of chance constraints. 
As in statistical estimation, the quality of our estimate is determined 
via confidence interval analysis. 

In contrast to existing approaches based on sampling, we
provide likelihood guarantees for the quality of the solutions found. In fact, we explicitly
state a confidence probability $\alpha$ that bounds the probability 
of exceeding a given error tolerance threshold $\vartheta$ in our estimation.
By properly choosing the estimation error $\vartheta$ and
the confidence probability $\alpha$ it is possible to generate compact sampled SCSPs
that can be effectively solved by existing solution methods. We also extended the reasoning
to SCOPs and demonstrated how to produce statistical upper and lower bounds for the value of the
optimal solution. 

We demonstrated our approach on a number of SCSPs and SCOPs: the static stochastic
knapsack problem, a stochastic multiprocessor scheduling problem, and a stochastic
lot-sizing problem. Our computational study demonstrates the effectiveness of our approach.

We conclude by briefly discussing a number of suggestions for future work.

\paragraph{Online stochastic optimization}
A promising direction is that of exploring synergies with online stochastic optimization \cite{bbcr01}. In particular, we suspect that our approach may be used to enhance the results in \cite{DBLP:conf/cpaior/HentenryckBV06,DBLP:conf/aips/MichelH04,DBLP:conf/aaai/BentH04,DBLP:conf/cp/BentKH05} by ensuring a better control of the solution quality obtained 
at each step of the online process.

\paragraph{Sampling strategies}
A key open issue is related to the fact that simple random sampling \cite{yates2002practice} is a relatively naive strategy for selecting samples. The use of more refined sampling strategies --- for instance a stratified sampling technique such as Latin Hypercube Sampling \cite{mbc79} --- may of course reduce the number of samples required to produce an ($\alpha$,$\vartheta$)-solution. Nevertheless, further research is required in order to clarify how stratified sampling can be effectively employed in this context.

\paragraph{Confidence intervals} The Clopper-Pearson interval is an exact interval since it is based directly on the binomial distribution rather than any approximation to the binomial distribution. This interval, however, can be conservative because of the discrete nature of the binomial distribution, as pointed out by Neyman \cite{Ney35}. For example, the true coverage rate of a 95\% Clopper-Pearson interval may be well above 95\%, depending on $n$ and $q$. Thus the interval may be wider than it needs to be to achieve 95\% confidence. In contrast, it is worth noting that other approximate confidence bounds may be narrower than their nominal confidence width, i.e., the ``normal approximation interval,'' also known as Wald confidence interval, the Wilson Interval, the Agresti-Coull Interval, etc, may in fact achieve a confidence level that is lower than the nominal one \cite{agresticloull98}. Future research may investigate the application of approximate intervals in the context of sample-based constraint solving. The performance of each of these approximate intervals have been thoroughly analysed in the existing body of literature. The advantage is that approximate intervals may lead to smaller sample sets and therefore to more compact sampled SCSPs.

\paragraph{Computational complexity} Finally, an interesting computational complexity questions remains open about the complexity of the standard multi-stage stochastic constraint programs.

\section*{Acknowledgments} 

The authors would like to thank the Associate Editor and the anonymous reviewers for commenting on our drafts and providing insightful remarks.

\appendix

\section{Filtering strategy for constraint expressions involving
expected values}
\label{sec:appendix_I}

We discuss a filtering strategy for handling constraint expressions involving
expected values in sampled SCSPs. This filtering strategy can be employed, 
in concert with the approach discussed in Section \ref{sec:scop}, to deal with 
the case in which the objective function is stochastic.
Consider a constraint $x=\text{E}[\langle \mathrm{exp} \rangle]$,
where $\text{E}[]$ denotes the expectation operator and $x$ is a real valued decision variable, whose domain
is stored as an interval with real valued upper and lower bounds.
Techniques for handling propagation and search involving
real valued decision variables are discussed in \cite{Benhamou06}.
A filtering algorithm that enforces bounds consistency on this
constraint is shown in Fig. \ref{sb_prop_alg_exp}.  
\begin{figure}[h!]
\vskip 1pc
  \begin{algorithm}[H]
    \SetKwInOut{Input}{input}
    \SetKwInOut{Output}{output}
    \Input{
    type; $\langle \text{exp} \rangle$; $\mathcal{T}$; $x$; $\alpha$.}
    \Output{Bound consistent $x$.}
    \BlankLine
    \Begin{
	$U \leftarrow \{\}$; $L \leftarrow \{\}$\;
        	\For{each  $p \in \Psi$}{
         		$U \leftarrow U \cup \mathrm{Sup}(\langle \text{exp} \rangle_{\downarrow{p}})$\;
         		$L \leftarrow L \cup \mathrm{Inf}(\langle \text{exp} \rangle_{\downarrow{p}})$\;
        	}
	$t\leftarrow \text{StudentT}(|\Psi|-1)$\;
	\If{\text{type}=$\mathcal{P}_{\text{lb}}$}{
			$\text{Sup}(x) \leftarrow \text{mean}(U)-\text{CDF}^{-1}_t(1-(1-\alpha)/2)\cdot\text{std}(U)/\sqrt{|\Psi|}$\;
			$\text{Inf}(x) \leftarrow \text{mean}(L)-\text{CDF}^{-1}_t(1-(1-\alpha)/2)\cdot\text{std}(L)/\sqrt{|\Psi|}$\;
	}
	\ElseIf{\text{type}=$\mathcal{P}_{\text{ub}}$}{
			$\text{Sup}(x) \leftarrow \text{mean}(U)+\text{CDF}^{-1}_t(1-(1-\alpha)/2)\cdot\text{std}(U)/\sqrt{|\Psi|}$\;
			$\text{Inf}(x) \leftarrow \text{mean}(L)+\text{CDF}^{-1}_t(1-(1-\alpha)/2)\cdot\text{std}(L)/\sqrt{|\Psi|}$\;
	}
    }
    \caption{Filtering Expected Values in sampled SCSPs \label{sb_prop_alg_exp}}
  \end{algorithm}
\end{figure}
It should be noted that the approach discussed in Section \ref{sec:scop} distinguishes two cases: the one in which our 
aim is to underestimates the true optimal profit (SCOP $\mathcal{P}_{\text{lb}}$) and that in which our aim is to overestimates the true optimal profit (SCOP $\mathcal{P}_{\text{ub}}$). The type of problem ($\mathcal{P}_{\text{lb}}$ or $\mathcal{P}_{\text{ub}}$) which the propagator belongs to must be specified as an input parameter ``type'' that influences propagation. 
The algorithm constructs two arrays: $U$ and $L$. $U$ lists, for each scenario, an upper bound for the expected value of $\langle \mathrm{exp} \rangle$,
$L$ lists, for each scenario, a lower bound for the expected value of $\langle \mathrm{exp} \rangle$. Then it exploits the Student's $t$ distribution with $|\Psi|-1$ degrees of freedom ($\text{StudentT}(|\Psi|-1)$)  to determine upper and lower confidence limits for the expected value of $\langle \mathrm{exp} \rangle$ at the prescribed confidence level $\alpha$. Note that $\text{CDF}^{-1}_t(\alpha)$ denotes the inverse cumulative distribution function of $t$; $\text{mean}(X)$ and $\text{std}(X)$ denote the mean and the standard deviation of the elements in $X$, respectively.
The algorithm operates by exploiting the structure $\Psi$ of the policy tree; therefore it takes implicitly into account the stage structure of the problem while computing the expected value of a given expression and it will correctly evaluate expected values both in a single or multi-stage case. Finally, it is worth remarking that this constraint is closely related to the Student's $t$ test constraint discussed in \cite{DBLP:conf/ecai/RossiPT14}.



\bibliographystyle{model1-num-names}
\bibliography{pub}







\end{document}